\pdfoutput=1
\documentclass[11pt,a4paper]{article}
\usepackage[authoryear,round]{natbib}
\usepackage{sjtu_iic_arxiv}
\usepackage{amsmath,amssymb,mathtools}
\usepackage{booktabs,multirow,array,longtable}
\usepackage{float}
\usepackage{placeins}
\usepackage{algorithmicx,algpseudocode}
\usepackage{graphicx}
\usepackage{xcolor}
\usepackage{microtype}
\usepackage{url}
\usepackage{hyperref}
\usepackage{etoc}
\usepackage{tikz}
\usepackage{enumitem}
\usetikzlibrary{positioning,calc}
\hypersetup{hidelinks}
\floatstyle{ruled}
\newfloat{algorithm}{tbp}{loa}
\floatname{algorithm}{Algorithm}
\newcolumntype{P}[1]{>{\raggedright\arraybackslash}p{#1}}
\newcolumntype{C}[1]{>{\centering\arraybackslash}p{#1}}
\definecolor{verdictblue}{HTML}{315A7D}
\definecolor{verdictgreen}{HTML}{2A9D8F}
\definecolor{verdictorange}{HTML}{E58F65}
\definecolor{verdictlight}{HTML}{F4F7FA}
\newtcolorbox{evidencecase}[1]{
  enhanced,
  width=\textwidth,
  colback=verdictlight,
  colframe=verdictblue,
  colbacktitle=verdictblue,
  coltitle=white,
  fonttitle=\small\bfseries,
  toptitle=0.6mm,
  bottomtitle=0.6mm,
  title={#1},
  boxrule=0.75pt,
  arc=2.2mm,
  outer arc=2.2mm,
  left=2.5mm,
  right=2.5mm,
  top=2.0mm,
  bottom=2.0mm,
  boxsep=0pt,
  before skip=8pt,
  after skip=10pt
}

\newcommand{\method}{\textsc{ORSMT}}

\newcommand{\refuted}{\textsc{Refuted}}

\title{\method: Unleashing the Potential of Satisfiability\\
Modulo Theories for Semantic Verification\\
and Repair in LLM-Based Optimization}
\setheadertext{%
  \begin{minipage}[c]{0.76\textwidth}%
  \begin{minipage}[c]{0.23\textwidth}\centering
    \includegraphics[trim=94bp 210bp 91bp 156bp,clip,width=\linewidth,height=0.95cm,keepaspectratio]{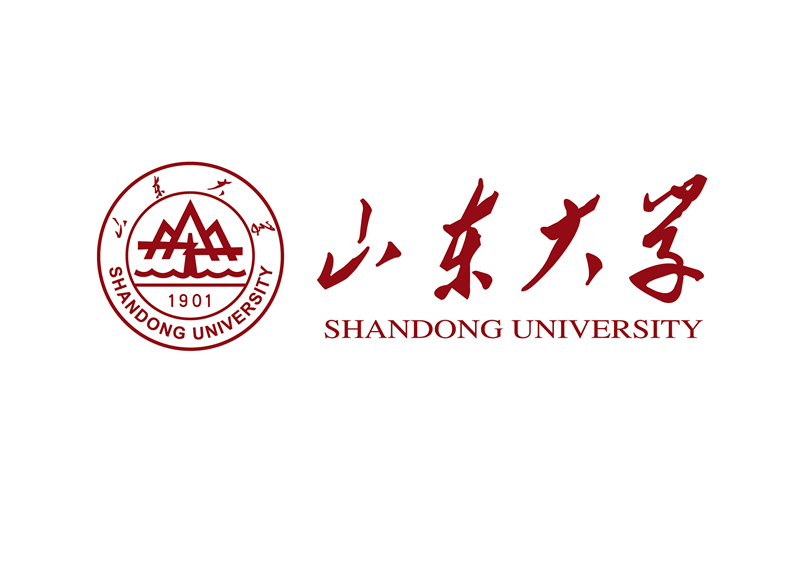}
  \end{minipage}\hfill
  \begin{minipage}[c]{0.23\textwidth}\centering
    \includegraphics[trim=76bp 67bp 73bp 73bp,clip,width=\linewidth,height=0.95cm,keepaspectratio]{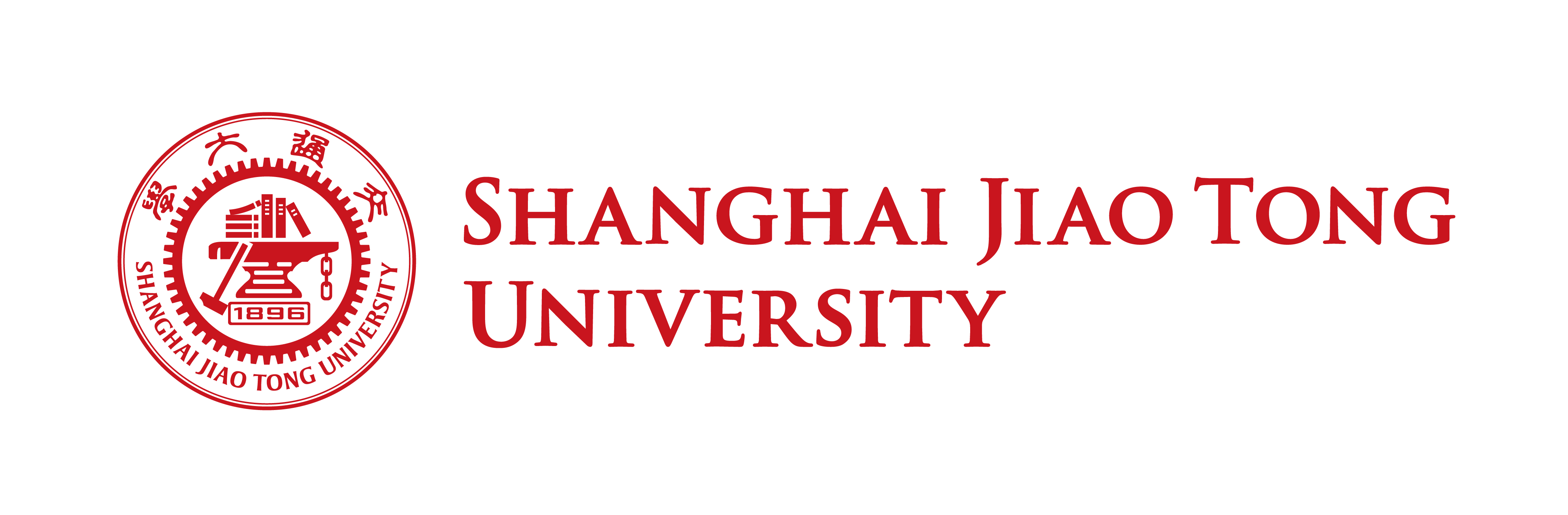}
  \end{minipage}\hfill
  \begin{minipage}[c]{0.27\textwidth}\centering
    \includegraphics[width=\linewidth,height=0.95cm,keepaspectratio]{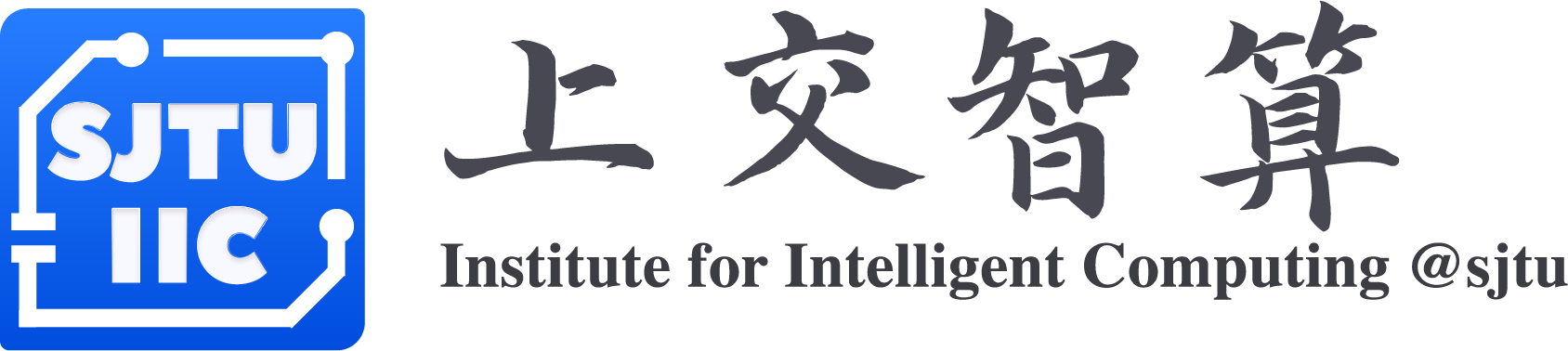}
  \end{minipage}\hfill
  \begin{minipage}[c]{0.22\textwidth}\centering
    \includegraphics[trim=174bp 174bp 300bp 174bp,clip,width=\linewidth,height=0.95cm,keepaspectratio]{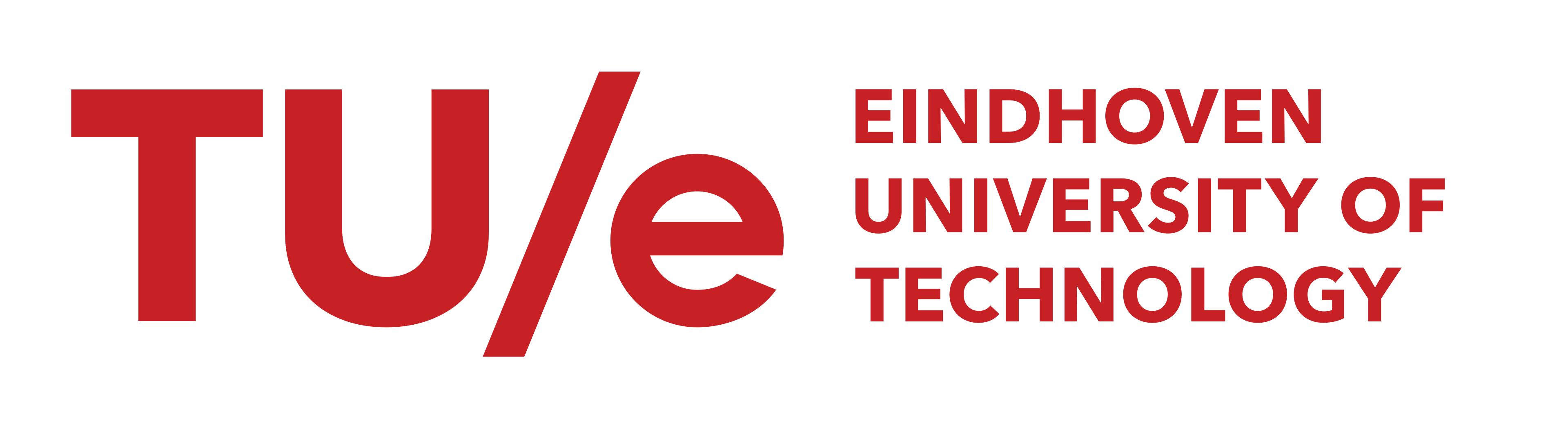}
  \end{minipage}%
  \end{minipage}%
}
\setheadertitle{ORSMT: Semantic Verification and Repair in LLM-Based Optimization}
\setversiontext{October 2026}
\githublink{https://github.com/ZRay1213/ORSMT}
\correspondingemail{%
  \href{mailto:202300190245@mail.sdu.edu.cn}{\nolinkurl{202300190245@mail.sdu.edu.cn}}\quad
  \href{mailto:wensong@email.sdu.edu.cn}{\nolinkurl{wensong@email.sdu.edu.cn}}\\[1pt]
  \texttt{\{\href{mailto:chenyuzhou@sjtu.edu.cn}{chenyuzhou},
  \href{mailto:linjianghao@sjtu.edu.cn}{linjianghao},
  \href{mailto:ddge@sjtu.edu.cn}{ddge}\}@sjtu.edu.cn}\quad
  \href{mailto:y.wu2@tue.nl}{\nolinkurl{y.wu2@tue.nl}}\\[2pt]
  \textsuperscript{\dag} Corresponding authors%
}
\author{%
  \mbox{Rui Zhang\textsuperscript{1},\hspace{0.3em}%
  Chenyu Zhou\textsuperscript{2},\hspace{0.3em}%
  Jianghao Lin\textsuperscript{2,\dag},\hspace{0.3em}%
  Yaoxin Wu\textsuperscript{3,\dag},\hspace{0.3em}%
  Dongdong Ge\textsuperscript{2},\hspace{0.3em}%
  Wen Song\textsuperscript{1}}\\[5pt]
  {\sjtuiicAffilFont
  \textsuperscript{1}Shandong University\\
  \textsuperscript{2}Shanghai Jiao Tong University\\
  \textsuperscript{3}Eindhoven University of Technology}
}
\hypersetup{
  pdftitle={ORSMT: Unleashing the Potential of Satisfiability Modulo Theories for Semantic Verification and Repair in LLM-Based Optimization},
  pdfauthor={Rui Zhang, Chenyu Zhou, Jianghao Lin, Yaoxin Wu, Dongdong Ge, Wen Song},
  pdfsubject={Semantic verification and repair for LLM-based optimization}
}

\begin{document}

\begin{abstract}
Large language models increasingly automate optimization modeling from natural language problems,
producing written mathematics, solver code and a returned report. Verifying one optimization problem
across these artifacts requires semantic comparison: unequal models can attain the same objective,
equivalent models can use different auxiliary variables, and mutually consistent artifacts can omit
the same source requirement. ORSMT provides a typed semantic interface and exact SMT relation queries.
It aligns observable decisions, projects auxiliary variables, proves or refutes feasible region and
objective relations, localizes refutations to implicated rows or domains and links source texts
when available. It then reconstructs one artifact from another, and each adopted candidate passes
the required relation checks.
Across 3,034 problems and six generator sizes, ORSMT improves the model--code match rate in all 72
benchmark and generator combinations before selection. Protected selection raises the aggregate
rate from 43.2\% to 88.2\% and the share that are both correct and matched from 25.2\% to 40.2\%. On the 1,905 refuted outputs,
certified voting under the shared fallback makes 50.9\% of outputs correct and matched, compared with 44.8\% for
ReLoop and 49.4\% for Opt-Verifier. A separate one-round source revision followed by exact
reconstruction restores a model--code match for 84.4\% of the refuted pairs.
\end{abstract}

\maketitle
\etocdepthtag.toc{body}

\section{Introduction}

Optimization models support logistics, production planning, and energy scheduling by optimizing decisions under operational constraints. Decision quality depends on how natural-language requirements are translated into mathematical formulations, executable solver models, and reported decisions. Recent large language model (LLM) systems generate these artifacts together, while evaluation commonly focuses on solver execution and the final objective value. Establishing \textbf{semantic agreement across these artifacts} therefore motivates verification and repair at the level of the optimization problem itself.

Classical criteria such as feasibility and answer correctness evaluate the returned decision without assessing the agreement across the underlying artifacts \citep{ramamonjison2022nl4opt,ahmadi2024optimus,huang2025orlm,jiang2025llmopt,yang2025optibench}. Different feasible regions can attain the same optimum, so a correct objective value does not identify the implemented problem. Conversely, equivalent models can differ in auxiliary variables, redundant constraints, and syntax. Agreement between written mathematics and executable code can also miss a source requirement from the natural language omitted by both. Output-based evaluation cannot distinguish these failures or establish whether a proposed repair restores the required relations.

Recent methods broaden feedback beyond the final value through behavioral tests, structural comparisons, model-based judgments, and counterexample searches \citep{zhang2026sacopt,li2026falsification,zhai2025equivamap,wang2025orgeval,lian2026reloop,liu2026optverifier}.
Although these methods provide useful repair signals, passing behavioral tests or receiving a favorable model-based judgment does not establish semantic agreement between the math formulation and its code implementation.
A program can produce the expected outputs on tested cases while admitting decisions excluded by the formulation, whereas structural comparisons can reject equivalent models that use different auxiliary constructions.
Verification must therefore examine whether the different artifacts admit the same feasible decisions and preserve their objective preferences.
It should also localize discrepancies to guide repair and verify that the repaired artifacts satisfy the required relations before adoption.

This paper proposes \textbf{ORSMT}, a satisfiability modulo theories (SMT)-based verification and repair framework for source requirements, written mathematics, instantiated solver models, and returned reports. 
ORSMT translates the written model and captures the solver model from code in a common exact representation, aligning observable decisions and projecting auxiliary variables. Within the affine rational core, SMT queries check admitted source requirements, feasible regions, objective relations, and reports for a fixed instance. A satisfiable violation query returns a localized counterexample; an unsatisfiable query certifies the corresponding relation. The evidence supports direct reconstruction between mathematics and code or guides revision, and every adopted candidate passes the required checks again. This keeps diagnosis, repair, and acceptance as distinct stages of one formally checked process. Our contributions are summarized as follows:

\textbf{(i)} We introduce ORSMT, a semantic verification framework for source requirements, mathematical formulations, solver codes, and returned reports. Its common representation supports exact SMT checks within the admitted affine rational core, providing certificates or localized counterexamples for decided relations.

\textbf{(ii)} We develop a verification and repair loop that combines bidirectional reconstruction with counterexample-guided revision. Each adopted candidate passes the required relation checks, connecting diagnosis, repair, and acceptance through the same verifier.

\textbf{(iii)} We evaluate ORSMT on 3,034 problems from 12 benchmarks and 6 generator sizes. Model--code match rates improve in all 72 benchmark--generator combinations before selection. Under protected selection, the aggregate match rate rises from 43.2\% to 88.2\%, and the share of outputs with both a correct answer and a model--code match rises from 25.2\% to 40.2\%. ORSMT further improves correct-and-matched outputs when integrated with existing OptiMUS-0.3 and OptiMind pipelines, supporting its use as a modular verification and repair layer.

\section{Related Work}

\textbf{Autoformulation and evaluation.} NL4Opt~\citep{ramamonjison2022nl4opt} studies formulation.
OptiMUS and OptiMUS-0.3 build modular formulation and execution
pipelines~\citep{ahmadi2024optimus,ahmaditeshnizi2024optimus03}, while OptiMind adds expert-guided training and inference~\citep{zhang2025optimind}. ORLM, LLMOPT and ORThought generate and revise
mathematical models, code and solver outputs~\citep{huang2025orlm,jiang2025llmopt,yang2025orthought}. MAMO,
OptiBench and OPT Engine broaden evaluation tasks~\citep{huang2025mamo,yang2025optibench,chen2026optengine}.
Feasibility and reference objectives assess final answers, but do not establish that the written
model, instantiated code and report describe the same optimization problem.

\textbf{Verification and repair.} ReLoop tests program behavior under parameter
changes~\citep{lian2026reloop}. Opt-Verifier judges structure and solutions~\citep{liu2026optverifier},
and SAC Opt uses semantic anchors for revision~\citep{zhang2026sacopt}. EquivaMap and ORGEval compare
variable mappings or graph structure~\citep{zhai2025equivamap,wang2025orgeval}. Falsification uses solver-based testing under controlled instance transformations, whereas FLARE verifies MILP reformulations using theorem proving~\citep{li2026falsification,robbins2026flare}. ORSMT repairs instance-specific mathematical, executable and report artifacts. Therefore, repair comparisons use ReLoop and Opt-Verifier.

\textbf{Formal semantic verification.} Denotational tests evaluate program behavior beyond textual
similarity~\citep{zhong2020semantic}. Z3 and cvc5 decide arithmetic formulas~\citep{demoura2008z3,
barbosa2022cvc5}, while neural symbolic methods translate proposed statements into formal
queries~\citep{jha2023neurosymbolic,pan2023logiclm,ye2023satlm,olausson2023linc}. ORSMT verifies source
requirements, feasible regions, objective preferences and reports in this way. Its evidence identifies
missing, additional and redundant restrictions and supports repair from either usable model.
\section{Problem Setting}\label{sec:problem}

An autoformulation pipeline maps a source statement \(s\) and instance data to written
mathematics, a solver program \(\widehat P\), and a report
\(R=(\hat x,\hat z,\operatorname{status})\), where \(\hat x\) and \(\hat z\) represent the returned solution and its objective value.
The written mathematics is lowered to the mathematical model \(A\).
The instantiated solver model is represented by
\(B=\mathcal H(\widehat P)\). A source adapter supplies formal requirements
\(C=\mathcal F_{\rm src}(s)\), each linked to its supporting
source text and variable roles, when that
adapter admits the source. Verification checks the available relations among \(C\), \(A\), \(B\), and \(R\).

Define a model as \(A=(X_A,U_A,\Phi_A,f_A,\delta_A)\), where \(X_A\) contains observable decisions,
\(U_A\) contains auxiliary variables, \(\Phi_A\) contains domains and constraints, \(f_A\) is the
objective, and \(\delta_A\in\{\min,\max\}\). The executable model \(B=(X_B,U_B,\Phi_B,f_B,\delta_B)\) is defined in the same way. A fixed correspondence aligns
the observable decisions of the two models. Their comparison projects out the auxiliary variables.
The projected feasible region and objective graph are
\begin{equation}
F_A^X=\{x:\exists u_A\ \Phi_A(x,u_A)\},\qquad
G_A^X=\{(x,z):\exists u_A\ \Phi_A(x,u_A)\land z=f_A(x,u_A)\}.
\label{eq:modelobjects}
\end{equation}
When \(G_A^X\) is single-valued, \(v_A(x)\) denotes the unique objective value associated with \(x\in F_A^X\).The relation \(x\succ_A y\)
denotes a strict improvement according to \(\delta_A\), and \(x\succeq_A y\) also permits equal
values. Therefore we define semantic model agreement as
\begin{equation}
A\equiv_{\rm sem}B\quad\Longleftrightarrow\quad(F_A^X,\succ_A)=(F_B^X,\succ_B).
\label{eq:modelagreement}
\end{equation}
This relation preserves feasible decisions and their preference order. Absolute objective values are
checked separately. Source satisfaction requires a feasible model that entails every applicable
requirement in a nonempty admitted \(C\). Untranslated source units retain an undecided source
verdict rather than a vacuous certificate. The report $R$ passes when its returned decision is feasible,
correctly valued and optimal. Instance parameters are substituted before comparison, and partial
reports use existential completion over unreported coordinates.

The task is to return matched artifacts and evidence for each required relation. A mathematical model
can supply a missing implementation, a captured solver model can supply missing mathematics, and a
source counterexample can expose an omission shared by both. Every decided relation is certified or
accompanied by a concrete refuting assignment.

\section{\method: Verification, Repair and Matching}\label{sec:method}

ORSMT first translates the  mathematical model $A$ and the executable model $B$ into a common exact representation, then uses SMT to check
their relations, and finally applies the same checks to every repaired candidate before adoption.

\begin{figure}[htbp]
\centering
\resizebox{\textwidth}{!}{\input{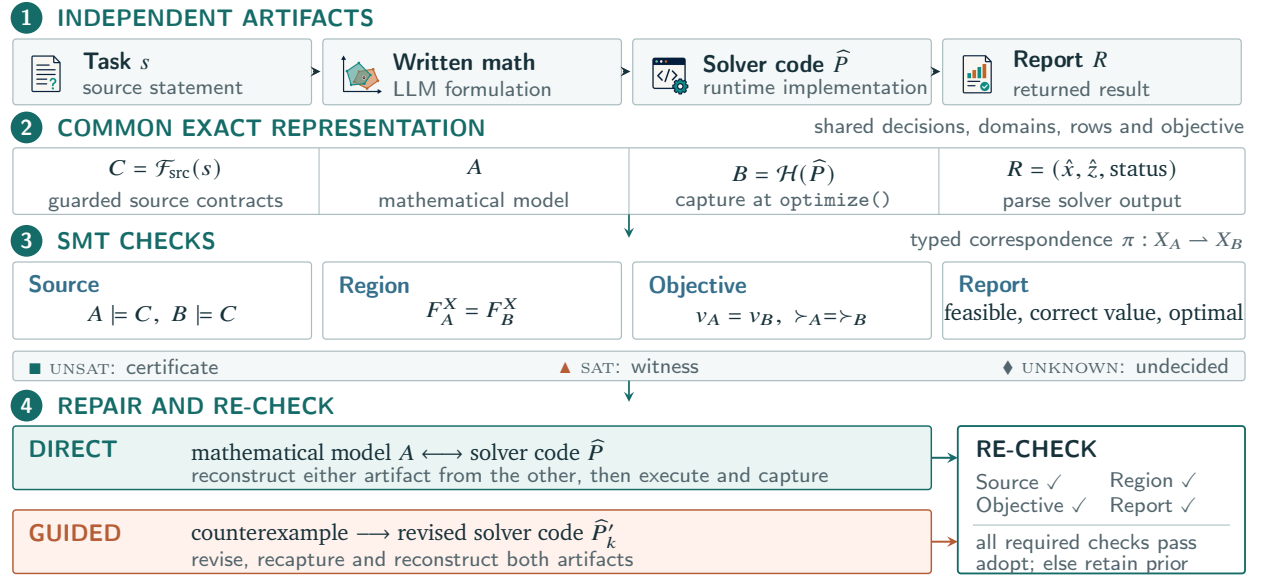}}
\caption{ORSMT translates source requirements, written mathematics, solver code and reports into a common exact representation. SMT checks either certify their relations or return a counterexample for direct reconstruction or revision. After repair, every relaton is checked again before adoption.}
\label{fig:completionflow}
\end{figure}

Figure~\ref{fig:completionflow} separates two models translation, SMT checking and the two repair routes.
The verifier supports affine models with rational coefficients over real, integer, and Boolean variables.
Appendix~\ref{app:algorithmsteps} gives the complete verification and repair procedure.

\subsection{Model representation and bidirectional construction}\label{subsec:lowering}

\textbf{Mathematical translation.} The translator \(\mathcal L\) lowers admitted mathematical formulation into the model
representation \(A\) defined in Section~\ref{sec:problem}. Structured declarations, solver objects, LP/MPS files provide
variables, contraints and objectives, 
while mathematical prose is transcribed together with links to the source text from which each element is
derived.
Finite decimal coefficients are represented as exact rationals. The resulting representation explicitly records variable bounds and domains, contraints, and contraint rows. Appendix~\ref{app:loweringdetails} presents the detailed input formats.

\textbf{Executable model and correspondence.} ORSMT aligns the observable
decisions $X_A$ and $X_B$ of A and B through a one-to-one correspondence based on normalized scalar
or index-tuple identifiers, variable roles, and types. Auxiliary variables $U_A$ and $U_B$ are not matched directly: they are projected out when comparing the model objects in
Eq.~\eqref{eq:modelobjects}. A contradiction, type mismatch, or unmatched observable decision prevents
a syntactic correspondence from being accepted and leaves the corresponding relation \textsc{Undecided}.
This separation allows semantically equivalent models to use different auxiliary constructions or
redundant constraints without requiring a one-to-one match of all solver variables.

\textbf{Target and normalization.}
ORSMT supports reconstruction in both directions. When the mathematical model \(A\) is available, ORSMT uses it as the primary construction target;
otherwise, it uses the captured executable model \(B=\mathcal H(\widehat P)\) as the target for reconstructing missing mathematical
artifacts. 
Before construction, ORSMT applies a normalization operator ($\mathcal N$) that preserves
the projected feasible region and objective graph relevant to the observable decisions. The
normalization may make implicit domains explicit and derive redundant constraints or auxiliary
definitions needed for compilation, while preserving the semantic model objects used for
verification. The resulting normalized representation is then used as the construction target.

\textbf{Mathematics to program.}
For the mathematical target \(T=A\), ORSMT computes
\begin{equation}
\widetilde A=T,\qquad
\widetilde P=\operatorname{Compile}(\mathcal N(T)),\qquad
\widetilde B=\mathcal H(\widetilde P),\qquad
\widetilde R=\operatorname{Solve}(\widetilde P).
\label{eq:construction}
\end{equation}
The generated solver code declares the normalized model and reports its
decision, objective value, and status. ORSMT then captures the instantiated executable model
($\widetilde B$) from the generated program and verifies its semantic agreement with
($\widetilde A$), together with the correctness of the returned report ($\widetilde R$).
Consequently, this direction can be used either to synthesize missing solver code from a mathematical
model or to repair existing code that fails to match the available formulation.

\textbf{Program to mathematics.} When the executable model is available but the mathematical formulation is missing or incomplete,
ORSMT sets (T=B). It exports the domains, bounds, constraints, objective, and optimization direction
of B into a mathematical representation ($\widetilde A$), and uses this representation as the
construction target. The resulting artifacts are compiled, executed, and captured as in Eq.~\eqref{eq:construction}, after which the same verification queries are applied to certify their
correspondence. Thus, either representation can serve as the construction source, while both
directions return to the same semantic verification procedure.

\textit{Proposition 1: Artifact Construction.}\label{prop:construction} If the normalization operator $\mathcal N$ preserves the projected feasible region and
objective graph of its input, and that $\mathcal H$ faithfully captures the model compiled from the
normalized representation,
Eq.~\eqref{eq:construction} constructs artifacts whose projected model semantics
agree with the construction target, up to the fixed correspondence of observable decisions. 
Appendix~\ref{app:formalguarantees} gives the full proof and the guarantees for
constraint reduction.


\subsection{Verification queries and diagnosis}\label{subsec:deciding2}

\textbf{Feasible regions.} ORSMT encodes domains and rows as exact arithmetic formulas and
computes the two directional violation queries
\begin{equation}
Q_{A\setminus B}=\operatorname{SAT}(F_A^X\land\neg F_B^X),\qquad
Q_{B\setminus A}=\operatorname{SAT}(F_B^X\land\neg F_A^X).
\label{eq:region}
\end{equation}
A witness in \(A\setminus B\) identifies a mathematical decision excluded by the implementation.
A witness in \(B\setminus A\) identifies an implemented decision excluded by the mathematics.
Unsatisfiability in both directions establishes equality of the projected
feasible regions which includes the case where both regions are empty. 
Therefore, model feasibility is checked separately before certifying semantic agreement.


\textbf{Verdict and diagnosis.}
Each verification query returns either a certificate or a counterexample. A satisfiable violation query refutes the corresponding relation and returns a concrete witness, whereas an unsatisfiable query certifies that the relation holds. After substituting instance parameters, ORSMT extracts the implicated constraints, domains, and linked source texts from the witness and, when applicable, the unsatisfiable core. Objective diagnoses additionally record the evaluated preference relation, while source diagnoses identify the violated requirement and its applicable scope. The resulting diagnostic payload is used for deterministic reconstruction or counterexample-guided language-model revision.

\textbf{Objective relations.} 
For decisions that are feasible in both models, ORSMT verifies whether the two objectives induce the same preference ordering. It decides the directional violation query:
\begin{equation}
\begin{split}
Q^{\rm pref}=\operatorname{SAT}\bigl(&x,y\in F_A^X\cap F_B^X\land\\
&[(x\succ_A y\land y\succeq_B x)\lor(x\succ_B y\land y\succeq_A x)]\bigr).
\end{split}
\label{eq:preference}
\end{equation}
A satisfying pair exposes opposite ordering or a strict improvement that becomes a tie. Feasible
models and unsatisfiable relation queries certify Eq.~\eqref{eq:modelagreement}. Absolute values use
\begin{equation}
Q^{\rm value}=\operatorname{SAT}\bigl(x\in F_A^X\cap F_B^X\land v_A(x)\ne v_B(x)\bigr).
\label{eq:objectivevalues}
\end{equation}
Unsatisfiability certifies equal values at common decisions, separately from preference agreement.


\textbf{Source requirements.}
Each source requirement in \(C\) is represented as a guarded contract
\((g_k,\psi_k,S_k)\), consisting of an applicability condition, a logical obligation, and its supporting source text. Variable identifiers are resolved through the typed source schema before a contract is admitted. Unresolved identifiers or untranslated source units remain undecided, and a source certificate is issued only when the admitted contract set is nonempty. For a feasible candidate, ORSMT decides
\begin{equation}
Q^{\rm source}_{A,k}=\operatorname{SAT}(\Phi_A\land g_k\land\neg\psi_k),\qquad
Q^{\rm source}_{B,k}=\operatorname{SAT}(\Phi_B\land g_k\land\neg\psi_k).
\label{eq:sourcequeries}
\end{equation}
Unsatisfiability certifies that the guarded requirement is enforced whenever its applicability condition holds, whereas a satisfying assignment returns a violating decision. A shared witness on both artifacts identifies a common omission and links it to \(S_k\). This source relation is independent of pairwise model equivalence: two artifacts may agree semantically while violating the same source requirement. 

\textbf{Missing, additional and redundant constraints.} Row entailment applies the same
violation test to individual mathematical and executable restrictions. A written row violated by an
implemented decision is missing from the implementation. An executable row violated by a written
decision adds a restriction. To test redundancy of a row \(c\), let \(D_A\) collect domains and
\(K_A\) collect constraint rows. ORSMT asks
\begin{equation}
Q^{\rm redundant}_{A,c}=
\operatorname{SAT}\!(D_A\land\bigwedge_{d\in K_A\setminus\{c\}}d\land\neg c).
\label{eq:redundancyquery}
\end{equation}
Unsatisfiability permits deletion of \(c\). Subsequent tests use the reduced model. A refuting
assignment is completed against the entire candidate before revision, and infeasible models receive
conflict diagnoses. Appendix~\ref{app:scale} details the reduction and scalability.

\textbf{Reported result.} Let \(O\subseteq X_B\) be the reported coordinates and \(\tau\geq0\) the
declared reporting tolerance. Let \(B_\tau(\hat x_O)\) as their tolerance box denote the corresponding tolerance box, with discrete
coordinates fixed to admissible integers; coordinates outside \(O\) are existentially completed.
For minimization, ORSMT computes
\begin{align}
Q^{\rm extend}_\tau&=\operatorname{SAT}(\Phi_B\land B_\tau(\hat x_O)),\nonumber\\
Q^{\rm bind}_\tau&=\operatorname{SAT}(\Phi_B\land B_\tau(\hat x_O)\land|f_B-\hat z|\leq\tau),\label{eq:solution}\\
Q^{\rm better}_\tau&=\operatorname{SAT}(\Phi_B\land f_B<\hat z-\tau).\label{eq:better}
\end{align}
Maximization reverses the last inequality. 
Unsatisfiability of the extension and binding queries certifies, respectively, that no feasible completion matches the reported decision and value, while satisfiability provides a completion witness.
An unsatisfiable better-decision query certifies that no feasible decision improves on \(\hat z\) beyond the declared tolerance.

\subsection{Certificates and scalable exactness}
\label{subsec:guarantees}

\textbf{Relation certificate.}
With fixed observable correspondence, feasible models, and single-valued objective graphs, unsatisfiability of both region-difference queries in Eq.~\eqref{eq:region} and the preference query in Eq.~\eqref{eq:preference} certifies \(A\equiv_{\rm sem}B\). A satisfying assignment to any query refutes the corresponding relation. The certificate covers all decisions admitted by the encoded models, while a refutation provides the separating decision or decision pair used for diagnosis. Appendix~\ref{app:formalguarantees} proves this claim and the corresponding report guarantees.

\textbf{Sound reduction.}
Using the domain formula \(D_A\) and active row set \(K_A\) defined above, consider a compiled obligation \(\psi\) and any \(L\subseteq K_A\):
\begin{equation}
\operatorname{UNSAT}(D_A\land\bigwedge_{c\in L}c\land\neg\psi)
\ \Longrightarrow\ \Phi_A\models\psi.
\label{eq:localrules}
\end{equation}
If the reduced query is unsatisfiable, the full violation query is also
unsatisfiable, certifying the obligation. If it is satisfiable, we test
whether its assignment can be extended to satisfy the omitted rows.
Only an extension that satisfies the full model and violates $\psi$
is reported as a counterexample; otherwise, the reduced query alone
does not refute the obligation.

\textbf{Exact component slicing.}
The bipartite row--variable incidence graph retains every connected component containing variables of \(\psi\). On a feasible base model, this slice is exact because disconnected components admit independent satisfying assignments. Smaller neighborhood slices retain the one-way guarantee in Eq.~\eqref{eq:localrules}. Appendix~\ref{app:formalguarantees} gives both proofs.

\subsection{Revision and Selection}

\textbf{Revision.}
Let $Y_0=(A,\widehat P,R)$ be the current artifacts. Verification evidence
$e$ records the failed relation, a counterexample when available, implicated
rows or domains, and linked source text. Given $s$ and $e$, the reviser
proposes a program $\widehat P'_k$. ORSMT captures its executable model,
reconstructs the artifacts through Eq.~(\ref{eq:construction}), and rechecks
the required source, region, objective, and report relations.

The evidence guides revision: a decision in the mathematical region but not
the executable region must be admitted by the implementation, while a decision
only in the executable region must be excluded. Objective and source revisions
address the violated preference and linked requirement, respectively. For
source revision, plan review checks the model against the source through the
computed plan; row review checks individual requirements to localize the
mismatch.

\textbf{Selection.}
Only candidates passing every check required by the active policy are
eligible. Direct selection takes the first eligible candidate in a fixed
order. Certified voting groups eligible candidates by verified reported
value and selects from the largest group if at least two agree. Protected
selection also preserves the original finite result when available. If no
candidate qualifies, ORSMT retains $Y_0$; otherwise, it adopts the selected
model, program, report, and verification evidence together. Reference
objectives and statuses are used only for evaluation. Appendix~B.4 gives
the complete procedure.

\textit{Proposition 2: Certified Adoption.}
Every adopted candidate satisfies the checks required by its selection
policy; an undecided required check prevents adoption. Otherwise, the
original artifacts remain unchanged. Voting additionally certifies agreement
of at least two eligible candidates on the reported value, and protected
selection preserves an available original finite result.

Checks follow their prerequisites: admission and feasibility precede semantic
comparison, region agreement precedes objective preference, and report
binding precedes optimality. An undecided prerequisite leaves its dependent
check undecided. 

Certificates are SMT-backed verdicts for the encoded relations, not
independently checked proof objects. They depend on translation, runtime
capture, observable correspondence, compiled source requirements, model
compilation, and the SMT backend. Unsupported translation, ambiguous
correspondence, timeout, or solver unknown yields \textsc{Undecided}.
Appendix~C.6 compares independent Z3 and cvc5 encodings.
\section{Experiments}\label{sec:experiments}

\textbf{Benchmarks and protocol.}\label{subsec:setup}\label{subsec:stage1efficiency}
The experiments evaluate SMT checking effectiveness and cost, direct reconstruction,
revision from counterexamples, repair baselines and host integration. A correct answer passes the
benchmark or host evaluator. A model--code match means that the written mathematical model and
captured executable model satisfy the exact semantic relation. An artifact match is the corresponding
reference-blind check for host integrations. An output is correct and matched when its answer is correct
and the applicable match holds. Problems decided is the percentage of the 3,034 assigned problems that return a decision; unavailable executable models count as unresolved. Failures repaired
is the percentage of original failures that are repaired. All checks pass means a model--code match,
feasible models and a verified report. Source requirements are evaluated separately.

The twelve benchmarks contain 3,034 problems. NL4Opt, NLP4LP, MAMO, OptiBench and Text2Opt provide
the main formulation collections~\citep{ramamonjison2022nl4opt,ahmadi2024optimus,huang2025mamo,
yang2025optibench,gao2026text2opt}. IndustryOR, LogiOR, ComplexOR, OptMATH and BWOR extend the
coverage~\citep{huang2025orlm,yang2025orthought,xiao2024chainofexperts,lu2025optmath,zhang2026bwor}.
For each of the 3,034 problems, one Qwen2.5-Instruct output at each of six sizes
(1.5B--72B) yields $3{,}034\times6=18{,}204$ assigned outputs~\citep{yang2024qwen25}.
References comprise 3,032 numeric objectives and two
infeasibility statuses (BWOR-037 and BWOR-054)~\citep{zhang2026bwor}. The controlled comparison fixes the 1,905
pre-repair refuted outputs, gives every feedback method the same Qwen2.5-14B reviser and permits
one temperature-zero revision with a 32,768-token context. A numeric
answer is correct when $|v-v^\star|\leq 10^{-6}\max(1,|v^\star|)$.
\begin{table}[htbp]
\caption{Complete assigned population: 3,034 problems for each of six generators, or 18,204 outputs.}
\label{tab:experimentalsetting}
\centering\footnotesize\setlength{\tabcolsep}{4pt}
\begin{tabular*}{\textwidth}{@{\extracolsep{\fill}}lrlr@{}}
\toprule
Benchmark & Problems & Benchmark & Problems\\
\midrule
MAMO easy & 652 & MAMO complex & 211\\
Text2Opt template & 630 & OptMATH & 132\\
OptiBench & 422 & IndustryOR & 98\\
Text2Opt unstructured & 248 & LogiOR & 76\\
NLP4LP & 242 & BWOR & 75\\
NL4Opt & 230 & ComplexOR & 18\\
\midrule
\multicolumn{3}{l}{\textbf{Total across 12 benchmarks}} & \textbf{3,034}\\
\bottomrule
\end{tabular*}
\end{table}

\textbf{SMT checking effectiveness.}\label{subsec:verifiercomparison}\label{subsec:refutecases}
The exact queries define 1,905 refuted and 7,867 certified-equivalent original pairs.
Table~\ref{tab:failure-surface-main} shows that code is unavailable nearly twice as often as
mathematics and contributes more unmatched restrictions.
\begin{table}[htbp]
\caption{Executable code adds a second failure surface over 3,034 problems and six generators (18,204 outputs). The final row uses 951 fully decided pairs; categories can overlap.}
\label{tab:failure-surface-main}
\centering\footnotesize
\setlength{\tabcolsep}{3pt}
\begin{tabular*}{0.88\textwidth}{@{\extracolsep{\fill}}lrr@{}}
\toprule
Audit & Written model & Executable model \\
\midrule
Artifact unavailable & 3,202 (17.6\%) & \textbf{5,980 (32.9\%)} \\
Primary availability failure & Malformed literal, 2,174 & \textbf{Runtime error, 3,723} \\
Redundant constraint on this side only & 188 & \textbf{239} \\
Constraint absent from the other side & 415 & \textbf{546} \\
\bottomrule
\end{tabular*}
\end{table}

\begin{table}[htbp]
\caption{Baseline flags on exact relation cohorts. Rates use available decisions in each cohort.}
\label{tab:verifiercoverage}
\centering\footnotesize\setlength{\tabcolsep}{4pt}
\begin{tabular*}{0.90\textwidth}{@{\extracolsep{\fill}}lrrrr@{}}
\toprule
Verifier & \shortstack{Refuted pairs\\flagged} & \shortstack{Refutation recall\\(\%)} &
\shortstack{Equivalent pairs\\flagged} & \shortstack{False positive rate\\(\%)}\\
\midrule
ReLoop execution & 181/1,897 & 9.5 & 29/7,859 & 0.4\\
Opt-Verifier structure & 1,656/1,905 & 86.9 & 4,576/7,867 & 58.2\\
\bottomrule
\end{tabular*}
\end{table}

Table~\ref{tab:verifiercoverage} compares baseline coverage on the same exact cohorts. ReLoop flags 181 available refutations and 29 equivalent pairs; Opt-Verifier flags 1,656 and 4,576, respectively.
Controlled semantic edits and independent Z3 and cvc5 encodings appear in
Appendices~\ref{app:criterionvalidity} and~\ref{app:executablescaling}. A separate source-covered
audit reports source-contract certificates, refutations and undecided pairs in
Appendix~\ref{app:sourcereview}.

\textbf{SMT efficiency and scaling.}\label{subsec:smtperformance}
Under a ten-second limit, automatic backend selection decides 94.6\% of the
3,034 problems with a 0.647-second median (Table~\ref{tab:smtscalingmain}).
On FrontierOR, exact component slicing yields paired speedups of
29.3$\times$ and 27.8$\times$ at 1,000--10,000 and 10,000--100,000 rows,
respectively; all 5,055 jointly decided claims agree.
Hardware details, query counts, and backend comparisons appear in
Appendix~\ref{app:executablescaling}.
\begin{table}[htbp]
\caption{ORSMT checking under a ten second limit on all 3,034 problems. The two settings run the same relation checks, and missing executable models are unresolved.}
\label{tab:smtscalingmain}
\centering
\footnotesize
\setlength{\tabcolsep}{7.0pt}
\begin{tabular*}{\textwidth}{@{\extracolsep{\fill}}lrrrr@{}}
\toprule
Setting & Problems & Decided & \shortstack{Problems\\decided (\%)} & \shortstack{Median checked\\time (s)} \\
\midrule
Complete audit, all problems & 3,034 & 2,860 & 94.3 & 0.397 \\
Automatic backend selection, all problems & 3,034 & 2,869 & 94.6 & 0.647 \\
\midrule
Automatic selection, no captured model & 134 & 0 & 0.0 & -- \\
Automatic selection, 1--4 constraints & 1,508 & 1,508 & 100.0 & 0.550 \\
Automatic selection, 5--24 constraints & 946 & 942 & 99.6 & 0.694 \\
Automatic selection, 25--99 constraints & 345 & 343 & 99.4 & 0.892 \\
Automatic selection, 100--999 constraints & 87 & 72 & 82.8 & 2.592 \\
Automatic selection, $\geq$1,000 constraints & 14 & 4 & 28.6 & 10.009 \\
\bottomrule
\end{tabular*}
\end{table}

\textbf{Direct reconstruction.}
Bidirectional reconstruction raises model--code matches from 7,867 to 15,724 without language model calls
by completing 6,590 programs from mathematics and 1,267 mathematical artifacts from programs.
The generator analysis retains every benchmark and model size before
protected selection.All 72 cells in Table~\ref{tab:completionmatrix} improve, with a median gain of 47.2 percentage points.
The match rate rises from 25.0\% to 93.4\% at 1.5B and from 14.0\% to 76.5\% at 3B; repaired 1.5B
outputs exceed unrepaired 72B outputs by 25.7 points.

\textbf{Revision from counterexamples and ablation.}\label{subsec:feedbackrepair}
Source revision raises outputs that are both correct and matched from 5,903 to 6,681 (Appendix Figure~\ref{fig:repairoutcomes}). On the 1,905 refuted outputs, the controlled ablation compares direct revision, plan review, row review and certified voting.
\begin{table}[htbp]
\caption{Model--code matches before and after repair over 3,034 problems per generator (\%), measured before protected selection.}
\label{tab:completionmatrix}
\centering\footnotesize\renewcommand{\arraystretch}{0.96}\setlength{\tabcolsep}{0.5pt}
\begin{tabular*}{\textwidth}{@{\extracolsep{\fill}}lrrrrrrrrrrrrr@{}}
\toprule
\multirow{2}{*}{Benchmark} & \multirow{2}{*}{Problems} & \multicolumn{2}{c}{1.5B} & \multicolumn{2}{c}{3B} & \multicolumn{2}{c}{7B} & \multicolumn{2}{c}{14B} & \multicolumn{2}{c}{32B} & \multicolumn{2}{c}{72B} \\
\cmidrule(lr){3-4}\cmidrule(lr){5-6}\cmidrule(lr){7-8}\cmidrule(lr){9-10}\cmidrule(lr){11-12}\cmidrule(lr){13-14}
 & & Before & After & Before & After & Before & After & Before & After & Before & After & Before & After \\
\midrule
MAMO easy & 652 & 24.5 & \textbf{99.7} & 22.9 & \textbf{92.6} & 77.3 & \textbf{100.0} & 74.4 & \textbf{99.7} & 96.0 & \textbf{100.0} & 97.7 & \textbf{100.0} \\
Text2Opt template & 630 & 22.1 & \textbf{92.7} & 0.3 & \textbf{52.2} & 12.1 & \textbf{41.7} & 15.2 & \textbf{52.5} & 12.9 & \textbf{61.6} & 17.3 & \textbf{63.7} \\
OptiBench & 422 & 37.9 & \textbf{99.8} & 27.5 & \textbf{90.0} & 59.0 & \textbf{99.1} & 44.5 & \textbf{95.5} & 80.8 & \textbf{100.0} & 93.1 & \textbf{100.0} \\
Text2Opt unstructured & 248 & 7.7 & \textbf{65.3} & 4.8 & \textbf{79.8} & 35.5 & \textbf{83.5} & 56.5 & \textbf{96.8} & 72.6 & \textbf{99.2} & 68.1 & \textbf{98.0} \\
NLP4LP & 242 & 37.6 & \textbf{99.6} & 25.2 & \textbf{87.2} & 55.0 & \textbf{99.2} & 41.3 & \textbf{94.2} & 84.3 & \textbf{100.0} & 93.4 & \textbf{100.0} \\
NL4Opt & 230 & 36.5 & \textbf{99.6} & 23.9 & \textbf{89.1} & 53.5 & \textbf{98.7} & 40.4 & \textbf{94.3} & 80.4 & \textbf{100.0} & 93.9 & \textbf{100.0} \\
MAMO complex & 211 & 17.5 & \textbf{96.7} & 2.4 & \textbf{78.2} & 28.0 & \textbf{85.3} & 44.5 & \textbf{90.0} & 54.5 & \textbf{89.1} & 57.3 & \textbf{85.8} \\
OptMATH & 132 & 3.8 & \textbf{86.4} & 2.3 & \textbf{28.0} & 9.1 & \textbf{59.8} & 20.5 & \textbf{71.2} & 34.1 & \textbf{89.4} & 27.3 & \textbf{87.9} \\
IndustryOR & 98 & 28.6 & \textbf{92.9} & 11.2 & \textbf{76.5} & 38.8 & \textbf{89.8} & 38.8 & \textbf{89.8} & 61.2 & \textbf{92.9} & 55.1 & \textbf{95.9} \\
LogiOR & 76 & 14.5 & \textbf{85.5} & 7.9 & \textbf{69.7} & 27.6 & \textbf{73.7} & 44.7 & \textbf{89.5} & 50.0 & \textbf{96.1} & 43.4 & \textbf{92.1} \\
BWOR & 75 & 30.7 & \textbf{78.7} & 6.7 & \textbf{66.7} & 48.0 & \textbf{85.3} & 50.7 & \textbf{94.7} & 58.7 & \textbf{97.3} & 64.0 & \textbf{94.7} \\
ComplexOR & 18 & 16.7 & \textbf{72.2} & 5.6 & \textbf{83.3} & 55.6 & \textbf{100.0} & 72.2 & \textbf{94.4} & 72.2 & \textbf{94.4} & 66.7 & \textbf{100.0} \\
\midrule
\textbf{All benchmarks} & \textbf{3,034} & \textbf{25.0} & \textbf{93.4} & \textbf{14.0} & \textbf{76.5} & \textbf{44.5} & \textbf{82.1} & \textbf{44.4} & \textbf{85.6} & \textbf{63.7} & \textbf{90.3} & \textbf{67.7} & \textbf{90.3} \\
\bottomrule
\end{tabular*}
\end{table}

A separate one-round comparison gives every method the same refuted outputs and 14B reviser. Source
revision followed by exact reconstruction restores an ORSMT model--code match for 1,607 outputs
(84.4\%); Appendix Table~\ref{tab:alignmentcost} reports the matched baseline outcomes and cost.

\textbf{Comparison with other repair methods.}
Certified voting makes 50.9\% of outputs correct and matched, 6.1 points above ReLoop and 1.4 above
Opt-Verifier. The clustered interval excludes zero only against ReLoop (Appendix
Table~\ref{tab:pairednew}); Tables~\ref{tab:nativepipelines} and~\ref{tab:newrepairquality}
report the full and paired pipelines, respectively.
ORSMT leads all three outcomes in Table~\ref{tab:nativepipelines}; Appendix~\ref{app:nativefull} gives the paired results.

\begin{table}[htbp]
\caption{Full-population outcomes over 3,034 problems and six generators (18,204 outputs). Outcome cells report count (rate); correct-answer counts use the 18,192 outputs with references.}
\label{tab:nativepipelines}
\centering\footnotesize\setlength{\tabcolsep}{1.0pt}
\begin{tabular*}{\textwidth}{@{\extracolsep{\fill}}lrrrrr@{}}
\toprule
Pipeline & \shortstack{Correct\\answer} & \shortstack{Model--code\\match} & \shortstack{Correct and\\matched} & \shortstack{Correct and matched\\/ correct answers} & \shortstack{Newly correct / original\\correct kept/total}\\
\midrule
Original & 4,931 (27.1) & 7,867 (43.2) & 4,588 (25.2) & 4,588/4,931 & 0 / 4,931/4,931\\
ReLoop & 5,760 (31.7) & 8,285 (45.5) & 4,744 (26.1) & 4,744/5,760 & 829 / 4,931/4,931\\
Opt-Verifier & 7,070 (38.9) & 4,612 (25.3) & 3,807 (20.9) & 3,807/7,070 & 2,139 / 4,420/4,931\\
\textbf{ORSMT} & \textbf{7,309 (40.2)} & \textbf{16,053 (88.2)} & \textbf{7,308 (40.2)} & \textbf{7,308/7,309} & \textbf{2,378 / 4,931/4,931}\\
\bottomrule
\end{tabular*}
\end{table}

\textbf{Integration with existing modeling systems.}\label{subsec:fulltransfer}
Both integrations use the same 3,034 assignments. OptiMUS-0.3 uses the exact model--code relation;
OptiMind uses an exact scalar extractor or reference-blind two-of-three auditors and rechecks every reconstruction.
\begin{table}[htbp]
\caption{Correct and matched outputs among 1,905 refuted outputs.}
\label{tab:newrepairquality}
\centering\footnotesize
\setlength{\tabcolsep}{7pt}
\begin{tabular*}{0.62\textwidth}{@{\extracolsep{\fill}}lr@{}}
\toprule
Policy & Correct and matched (\%)\\
\midrule
ReLoop & 44.8\\
Opt-Verifier & 49.4\\
ORSMT direct revision & 46.4\\
ORSMT with plan review & 47.5\\
ORSMT with row review & 49.6\\
\textbf{ORSMT certified voting} & \textbf{50.9}\\
\bottomrule
\end{tabular*}
\end{table}

\begin{table}[htbp]
\caption{Host integration over 3,034 problems. Cells report count (rate).}
\label{tab:plugincomparison}
\centering
\footnotesize
\setlength{\tabcolsep}{3.0pt}
\begin{tabular*}{\textwidth}{@{\extracolsep{\fill}}llrrr@{}}
\toprule
System & Pipeline & \shortstack{Correct\\answer} &
\shortstack{Artifact\\match} & \shortstack{Correct and\\matched}\\
\midrule
OptiMUS-0.3 & Original & 647 (21.3) & 1,113 (36.7) & 630 (20.8)\\
 & Original + ORSMT & \textbf{1,199 (39.5)} & \textbf{1,938 (63.9)} & \textbf{1,199 (39.5)}\\
\midrule
OptiMind & Original & 1,885 (62.1) & 2,254 (74.3) & 1,678 (55.3)\\
 & Original + ORSMT & \textbf{1,904 (62.8)} & \textbf{2,937 (96.8)} & \textbf{1,901 (62.7)}\\
\bottomrule
\end{tabular*}
\end{table}

ORSMT adds 552 correct answers and 569 correct-and-matched outputs for OptiMUS-0.3,
versus 19 and 223 for OptiMind, respectively. OptiMind also gains 683 artifact matches,
showing that the integration chiefly improves artifact agreement when the host's answers
are already strong. Appendix~\ref{app:optimindaudit} gives routes, original-failure
repair rates and cost.

\section{Conclusion}

ORSMT closes the verification--repair loop by combining exact certificates and concrete counterexamples with bidirectional reconstruction and rechecking every adopted artifact. Across 3,034 problems and six generator sizes (18,204 outputs), it raises the model--code match rate from 43.2\% to 88.2\% and the share of outputs that are both correct and matched from 25.2\% to 40.2\%. Among 1,905 refuted outputs, certified voting yields 50.9\% correct-and-matched outputs, while one source revision followed by reconstruction restores model--code agreement for 84.4\% of refuted pairs. The analysis also reveals systematic asymmetries in generated artifacts: code is unavailable nearly twice as often as mathematics and is more likely to introduce additional constraints than to omit existing ones. ORSMT decides 2,869 of 3,034 problems, remains sub-second through 99 constraints, and achieves a 29.3$\times$ speedup from component slicing. Its gains on both host systems further support its use as a modular verification and repair layer for LLM-based optimization pipelines.

\bibliography{references}

@article{ahmaditeshnizi2024optimus03,
  title = {{OptiMUS-0.3}: Using Large Language Models to Model and Solve Optimization Problems at Scale},
  author = {AhmadiTeshnizi, Ali and Gao, Wenzhi and Brunborg, Herman and Talaei, Shayan and Lawless, Connor and Udell, Madeleine},
  journal = {arXiv preprint arXiv:2407.19633},
  year = {2024},
  eprint = {2407.19633},
  archivePrefix = {arXiv}
}

@article{zhang2025optimind,
  title = {{OptiMind}: Teaching {LLM}s to Think Like Optimization Experts},
  author = {Zhang, Xinzhi and Chen, Zeyi and Zope, Humishka and Barbalho, Hugo and Mellou, Konstantina and Molinaro, Marco and Kulkarni, Janardhan and Menache, Ishai and Li, Sirui},
  journal = {arXiv preprint arXiv:2509.22979},
  year = {2025},
  eprint = {2509.22979},
  archivePrefix = {arXiv}
}

@article{yang2024qwen25,
  title   = {{Qwen2.5} Technical Report},
  author  = {Yang, An and Yang, Baosong and Zhang, Beichen and Hui, Binyuan and Zheng, Bo and Yu, Bowen and Li, Chengyuan and Liu, Dayiheng and Huang, Fei and Wei, Haoran and others},
  journal = {arXiv preprint arXiv:2412.15115},
  year    = {2024},
  eprint  = {2412.15115},
  archivePrefix = {arXiv}
}

@inproceedings{xiao2024chainofexperts,
  title     = {Chain-of-Experts: When {LLM}s Meet Complex Operations Research Problems},
  author    = {Xiao, Ziyang and Zhang, Dongxiang and Wu, Yangjun and Xu, Lilin and Wang, Yuan Jessica and Han, Xiongwei and Fu, Xiaojin and Zhong, Tao and Zeng, Jia and Song, Mingli and Chen, Gang},
  booktitle = {International Conference on Learning Representations},
  year      = {2024},
  url       = {https://openreview.net/forum?id=HobyL1B9CZ}
}

@inproceedings{lu2025optmath,
  title     = {{OptMATH}: A Scalable Bidirectional Data Synthesis Framework for Optimization Modeling},
  author    = {Lu, Hongliang and Xie, Zhonglin and Wu, Yaoyu and Ren, Can and Chen, Yuxuan and Wen, Zaiwen},
  booktitle = {Proceedings of the 42nd International Conference on Machine Learning},
  year      = {2025},
  eprint    = {2502.11102},
  archivePrefix = {arXiv}
}

@inproceedings{ramamonjison2022nl4opt,
  title     = {Augmenting Operations Research with Auto-Formulation of Optimization Models from Problem Descriptions},
  author    = {Ramamonjison, Rindra and Li, Haley and Yu, Timothy T. and He, Shiqi and Rengan, Vishnu and Banitalebi-Dehkordi, Amin and Zhou, Zirui and Zhang, Yong},
  booktitle = {Proceedings of the 2022 Conference on Empirical Methods in Natural Language Processing: Industry Track},
  pages     = {29--62},
  year      = {2022},
  publisher = {Association for Computational Linguistics},
  doi       = {10.18653/v1/2022.emnlp-industry.4},
  eprint    = {2209.15565},
  archivePrefix = {arXiv}
}

@inproceedings{ahmadi2024optimus,
  title     = {{OptiMUS}: Scalable Optimization Modeling with {(MI)LP} Solvers and Large Language Models},
  author    = {AhmadiTeshnizi, Ali and Gao, Wenzhi and Udell, Madeleine},
  booktitle = {Proceedings of the 41st International Conference on Machine Learning},
  year      = {2024},
  eprint    = {2310.06116},
  archivePrefix = {arXiv}
}

@article{huang2025orlm,
  title   = {{ORLM}: A Customizable Framework in Training Large Models for Automated Optimization Modeling},
  author  = {Huang, Chenyu and Tang, Zhengyang and Hu, Shixi and Jiang, Ruoqing and Zheng, Xin and Ge, Dongdong and Wang, Benyou and Wang, Zizhuo},
  journal = {Operations Research},
  volume  = {73},
  number  = {6},
  pages   = {2986--3009},
  year    = {2025},
  doi     = {10.1287/opre.2024.1233},
  eprint  = {2405.17743},
  archivePrefix = {arXiv}
}

@inproceedings{jiang2025llmopt,
  title     = {{LLMOPT}: Learning to Define and Solve General Optimization Problems from Scratch},
  author    = {Jiang, Caigao and Shu, Xiang and Qian, Hong and Lu, Xingyu and Zhou, Jun and Zhou, Aimin and Yu, Yang},
  booktitle = {International Conference on Learning Representations},
  year      = {2025},
  eprint    = {2410.13213},
  archivePrefix = {arXiv}
}

@inproceedings{yang2025optibench,
  title     = {{OptiBench} Meets {ReSocratic}: Measure and Improve {LLM}s for Optimization Modeling},
  author    = {Yang, Zhicheng and Wang, Yiwei and Huang, Yinya and Guo, Zhijiang and Shi, Wei and Han, Xiongwei and Feng, Liang and Song, Linqi and Liang, Xiaodan and Tang, Jing},
  booktitle = {International Conference on Learning Representations},
  year      = {2025},
  eprint    = {2407.09887},
  archivePrefix = {arXiv}
}

@inproceedings{huang2025mamo,
  title     = {{LLM}s for Mathematical Modeling: Towards Bridging the Gap between Natural and Mathematical Languages},
  author    = {Huang, Xuhan and Shen, Qingning and Hu, Yan and Gao, Anningzhe and Wang, Benyou},
  booktitle = {Findings of the Association for Computational Linguistics: NAACL 2025},
  pages     = {2678--2710},
  year      = {2025},
  publisher = {Association for Computational Linguistics},
  doi       = {10.18653/v1/2025.findings-naacl.146},
  eprint    = {2405.13144},
  archivePrefix = {arXiv}
}

@inproceedings{zhai2025equivamap,
  title     = {{EquivaMap}: Leveraging {LLM}s for Automatic Equivalence Checking of Optimization Formulations},
  author    = {Zhai, Haotian and Lawless, Connor and Vitercik, Ellen and Leqi, Liu},
  booktitle = {Proceedings of the 42nd International Conference on Machine Learning},
  pages     = {74288--74305},
  volume    = {267},
  series    = {Proceedings of Machine Learning Research},
  publisher = {PMLR},
  year      = {2025},
  url       = {https://proceedings.mlr.press/v267/zhai25a.html},
  eprint    = {2502.14760},
  archivePrefix = {arXiv}
}

@article{robbins2026flare,
  title   = {{{FLARE}}: Verifying {{MILP}} Reformulations with {{LLM}}-Based Theorem Proving},
  author  = {Robbins, Henry and Lawless, Connor and Udell, Madeleine and Vitercik, Ellen},
  journal = {arXiv preprint arXiv:2608.25220},
  year    = {2026},
  eprint  = {2608.25220},
  archivePrefix = {arXiv}
}

@article{wang2025orgeval,
  title   = {{ORGEval}: Graph-Theoretic Evaluation of {LLM}s in Optimization Modeling},
  author  = {Wang, Zhuohan and Zhu, Ziwei and Li, Ziniu and Chen, Congliang and Han, Yizhou and Lin, Yufeng and Lin, Zhihang and Gu, Angyang and Hu, Xinglin and Sun, Ruoyu and Ding, Tian},
  journal = {arXiv preprint arXiv:2510.27610},
  year    = {2025},
  eprint  = {2510.27610},
  archivePrefix = {arXiv}
}

@article{jha2023neurosymbolic,
  title   = {Neuro Symbolic Reasoning for Planning: Counterexample Guided Inductive Synthesis using Large Language Models and Satisfiability Solving},
  author  = {Jha, Sumit Kumar and Jha, Susmit and Lincoln, Patrick and Bastian, Nathaniel D. and Velasquez, Alvaro and Ewetz, Rickard and Neema, Sandeep},
  journal = {arXiv preprint arXiv:2309.16436},
  year    = {2023},
  eprint  = {2309.16436},
  archivePrefix = {arXiv}
}

@inproceedings{pan2023logiclm,
  title     = {{Logic-LM}: Empowering Large Language Models with Symbolic Solvers for Faithful Logical Reasoning},
  author    = {Pan, Liangming and Albalak, Alon and Wang, Xinyi and Wang, William Yang},
  booktitle = {Findings of the Association for Computational Linguistics: EMNLP 2023},
  pages     = {3806--3824},
  year      = {2023},
  publisher = {Association for Computational Linguistics},
  doi       = {10.18653/v1/2023.findings-emnlp.248},
  eprint    = {2305.12295},
  archivePrefix = {arXiv}
}

@inproceedings{ye2023satlm,
  title     = {{SATLM}: Satisfiability-Aided Language Models Using Declarative Prompting},
  author    = {Ye, Xi and Chen, Qiaochu and Dillig, Isil and Durrett, Greg},
  booktitle = {Advances in Neural Information Processing Systems},
  volume    = {36},
  pages     = {45548--45580},
  year      = {2023},
  eprint    = {2305.09656},
  archivePrefix = {arXiv}
}

@inproceedings{olausson2023linc,
  title     = {{LINC}: A Neurosymbolic Approach for Logical Reasoning by Combining Language Models with First-Order Logic Provers},
  author    = {Olausson, Theo X. and Gu, Alex and Lipkin, Benjamin and Zhang, Cedegao E. and Solar-Lezama, Armando and Tenenbaum, Joshua B. and Levy, Roger},
  booktitle = {Proceedings of the 2023 Conference on Empirical Methods in Natural Language Processing},
  pages     = {5153--5176},
  year      = {2023},
  publisher = {Association for Computational Linguistics},
  doi       = {10.18653/v1/2023.emnlp-main.313},
  eprint    = {2310.15164},
  archivePrefix = {arXiv}
}

@inproceedings{demoura2008z3,
  title     = {{Z3}: An Efficient {SMT} Solver},
  author    = {de Moura, Leonardo and Bj{\o}rner, Nikolaj},
  booktitle = {Tools and Algorithms for the Construction and Analysis of Systems},
  pages     = {337--340},
  year      = {2008},
  publisher = {Springer},
  doi       = {10.1007/978-3-540-78800-3_24}
}

@inproceedings{barbosa2022cvc5,
  title     = {{cvc5}: A Versatile and Industrial-Strength {SMT} Solver},
  author    = {Barbosa, Haniel and Barrett, Clark and Brain, Martin and Kremer, Gereon and Lachnitt, Hanna and Mann, Makai and Mohamed, Abdalrhman and Mohamed, Mudathir and Niemetz, Aina and N{\"o}tzli, Andres and Ozdemir, Alex and Preiner, Mathias and Reynolds, Andrew and Sheng, Ying and Tinelli, Cesare and Zohar, Yoni},
  booktitle = {Tools and Algorithms for the Construction and Analysis of Systems},
  pages     = {415--442},
  year      = {2022},
  publisher = {Springer},
  doi       = {10.1007/978-3-030-99524-9_24}
}

@article{kong2026frontieror,
  title   = {{FrontierOR}: Benchmarking {LLM}s' Capacity for Efficient Algorithm Design in Large-Scale Optimization},
  author  = {Kong, Minwei and Jiang, Chonghe and Qu, Ao and Ouyang, Wenbin and Zeng, Zhaoming and Guo, Xiaotong and Li, Zhekai and Li, Junyi and Fan, Yi and Zheng, Xinshou and others},
  journal = {arXiv preprint arXiv:2605.25246},
  year    = {2026},
  eprint  = {2605.25246},
  archivePrefix = {arXiv}
}

@article{li2026falsification,
  title   = {Falsification-Based Verification of {LLM}-Generated Optimization Models: Sound Test Batteries and Their Detection Limits},
  author  = {Li, Haifeng and Hai, Mo},
  journal = {arXiv preprint arXiv:2607.16646},
  year    = {2026},
  eprint  = {2607.16646},
  archivePrefix = {arXiv}
}

@inproceedings{liu2026optverifier,
  title   = {{Opt-Verifier}: Unleashing the Power of {LLM}s for Optimization Modeling via Dual-Side Verification},
  author  = {Liu, Haoyang and Wang, Jie and Niu, Boxuan and Han, Xiongwei and Xu, Yian and Ye, Mingxuan and Geng, Zijie and Zhu, Fangzhou and Zhong, Tao and Yuan, Mingxuan and Hao, Jianye},
  booktitle = {Proceedings of the 43rd International Conference on Machine Learning},
  year    = {2026},
  eprint  = {2605.29556},
  archivePrefix = {arXiv}
}

@article{lian2026reloop,
  title   = {{ReLoop}: Structured Modeling and Behavioral Verification for Reliable {LLM}-Based Optimization},
  author  = {Lian, Junbo Jacob and Sun, Yujun and Chen, Huiling and Zhang, Chaoyu and Qin, Hanzhang and Teo, Chung-Piaw},
  journal = {arXiv preprint arXiv:2602.15983},
  year    = {2026},
  eprint  = {2602.15983},
  archivePrefix = {arXiv}
}

@inproceedings{zhang2026sacopt,
  title     = {{SAC-Opt}: Semantic Anchors for Iterative Correction in Optimization Modeling},
  author    = {Zhang, Yansen and Kang, Qingcan and Chen, Yujie and Wang, Yufei and Han, Xiongwei and Zhong, Tao and Yuan, Mingxuan and Ma, Chen},
  booktitle = {Proceedings of the 43rd International Conference on Machine Learning},
  year      = {2026},
  eprint    = {2510.05115},
  archivePrefix = {arXiv}
}

@misc{zhang2026bwor,
  title        = {{BWOR}: A Bilingual Benchmark for Modeling and Solving Operations Research Problems from Natural Language},
  author       = {Zhang, Bowen},
  year         = {2026},
  publisher    = {Zenodo},
  doi          = {10.5281/zenodo.20120692},
  url          = {https://doi.org/10.5281/zenodo.20120692},
  howpublished = {Dataset}
}

@article{yang2025orthought,
  title   = {{ORThought}: Benchmarking and Automating Logistics Optimization Modeling via Structured {LLM} Reasoning},
  author  = {Yang, Beinuo and Zhou, Qishen and Li, Junyi and Su, Chenxing and Angeloudis, Panagiotis and Hu, Simon},
  journal = {Artificial Intelligence for Transportation},
  volume  = {6},
  pages   = {100059},
  year    = {2026},
  doi     = {10.1016/j.ait.2026.100059},
  eprint  = {2508.14410},
  archivePrefix = {arXiv}
}

@inproceedings{gao2026text2opt,
  title   = {Models Can Model, But Can't Bind: Structured Grounding in Text-to-Optimization},
  author  = {Gao, Zhiqi and Ge, Albert and Berenbeim, Alexander and Bastian, Nathaniel D. and Sala, Frederic},
  booktitle = {Conference on Language Modeling},
  year    = {2026},
  eprint  = {2605.21751},
  archivePrefix = {arXiv}
}

@inproceedings{chen2026optengine,
  title     = {{OPT-Engine}: Benchmarking the Limits of {LLM}s in Optimization Modeling via Complexity Scaling},
  author    = {Chen, Yitian and Cheng, Cheng and Sun, Yinan and Ling, Zi and Ge, Dongdong},
  booktitle = {Proceedings of the 43rd International Conference on Machine Learning},
  year      = {2026},
  eprint    = {2601.19924},
  archivePrefix = {arXiv}
}

@inproceedings{zhong2020semantic,
  title     = {Semantic Evaluation for Text-to-SQL with Distilled Test Suites},
  author    = {Zhong, Ruiqi and Yu, Tao and Klein, Dan},
  booktitle = {Proceedings of the 2020 Conference on Empirical Methods in Natural Language Processing},
  pages     = {396--411},
  year      = {2020},
  publisher = {Association for Computational Linguistics},
  doi       = {10.18653/v1/2020.emnlp-main.29}
}
\bibliographystyle{plainnat}

\clearpage
\appendix
\etocdepthtag.toc{appendix}
\begingroup
\etocsettagdepth{body}{none}
\etocsettagdepth{appendix}{subsection}
\etocobeydepthtags
\etocsetstyle{section}
  {\par\smallskip}{}
  {\noindent{\bfseries\etocnumber\quad\etocname\hfill\etocpage}\par}{}
\etocsetstyle{subsection}
  {}{}
  {\noindent\hspace{1em}\etocnumber\quad\etocname\hfill\etocpage\par}{}
\etocsettocstyle{}{}
\small
\noindent\textbf{Contents of the Appendix}\par
\smallskip
\tableofcontents
\endgroup
The appendix follows the three experimental levels used by the paper. After fixing the populations,
representations, scoring rules and formal procedures, Section~\ref{app:scale} reports SMT checking
coverage, unresolved cases, error distributions and running time. Section~\ref{app:repair} then
separates feedback quality from repair effects before reporting full before-and-after results and
transfer to existing modeling systems. The repair, source-review and host-audit prompts appear beside
the experiments that use them. The final sections test the evaluation criterion and give worked cases
with exact witnesses.

\section{Experimental Protocol and Artifact Representations}\label{app:setup}\label{app:repairmetrics}
\subsection{Populations, metrics and scoring}
The main experiment pairs all twelve benchmarks with six Qwen2.5-Instruct generators. Each
problem retains its identifier, benchmark, generator, original artifacts, verifier version and
adopted output. The repair cohort consists of the originally refuted pairs, and source problem
resampling keeps outputs of the same problem together. Table~\ref{tab:supppopulations} distinguishes
the populations used in the main results and supplementary controls.

\begin{table}[htbp]
\caption{Assigned populations and units used in the experiments.}
\label{tab:supppopulations}
\centering\footnotesize\setlength{\tabcolsep}{3pt}
\begin{tabular*}{\textwidth}{@{\extracolsep{\fill}}llr@{}}
\toprule
Experiment & Unit & Assigned units \\
\midrule
Main generation and repair & Problem and generator & $3{,}034\times6=18{,}204$ \\
Paired feedback intervention & Originally refuted output & 1,905 \\
Paired constraint differences & Output with four decided comparisons & 951 \\
Controlled CPU verification & Unique problem & 3,034 \\
Canonical requirement queries & Model obligation & 136,100 \\
Canonical redundancy & Model row & 340,289 \\
Labeled edit controls & Edited or rewritten pair & 1,360 \\
Reference edit stress test & Edited pair & 12,831 \\
\bottomrule
\end{tabular*}
\end{table}

We report the four outcomes defined in Section~\ref{subsec:setup}: correct answer,
model--code match, artifact match and correct and matched. A numeric answer is correct when
\(\lvert z-z^\star\rvert\leq10^{-6}\max\{1,\lvert z^\star\rvert\}\). An output passes all checks when
the models match, both are feasible and the report is verified. Empty feasible regions
are counted separately in reconstruction outcomes. Every paired comparison preserves its assigned
identifiers, including unavailable artifacts, unresolved checks and execution failures.

The original verification panel uses exact decimal arithmetic for report discrepancies. Repair
scoring uses the declared tolerance for printed decisions and reference values.

\subsection{Artifact representations and capture}\label{app:loweringdetails}
Accepted declarations contain finitely many uniquely named real, integer or Boolean variables,
optional rational bounds, linear rows with relation \(\leq\), \(\geq\) or equality, and an objective
with its optimization direction. Structured declarations, solver objects, LP and MPS files supply
these fields directly. Mathematical prose supplies a proposed transcription whose source text is
retained. Malformed declarations, unsupported operators and undeclared references return a typed
translation failure.

The historical panel represents each variable with a real arithmetic symbol and adds integrality
predicates for integer domains. Boolean domains intersect integrality with the bounds zero and one.
An absent bound remains unbounded, and each finite decimal denotes its exact rational value.
Its accepted declarations require a nonempty objective term list and nonempty terms in each row.
The row list itself may be empty. Region and preference comparisons use zero tolerance.

Runtime capture intercepts the optimization call, flushes pending modifications and records the
instantiated model after loops and instance expansion. Later solver calls produce separate
snapshots. Returned text, explicit objective and decisions, solver attributes and execution
exceptions occupy separate fields. The frozen model comparison aligns all declared variables by
unique normalized names, with no hidden auxiliaries. Missing or ambiguous correspondences retain
an unresolved state.

\textbf{Source requirements.} Each obligation records its supporting source text, quantified scope,
variable roles, exact coefficients, relation and applicability. Released reference formulations
supply formal obligations when their adapters admit them. Open natural language supplies proposed
obligations for review. Reference conventions such as symmetry constraints retain their provenance.
Published objective values and status labels remain separate and enter only scoring. When an
adapter admits a released formulation as a source contract, the resulting statement only checks that
contract. It does not independently recover the author's intent.

\textbf{Exact reconstruction.} The compiler emits domains, rows, objective and direction from the
selected model. A bounded product \(y=bx\), where \(b\) is Boolean and \(L\leq x\leq U\), has the exact
reformulation
\[
Lb\leq y\leq Ub,\qquad x-U(1-b)\leq y\leq x-L(1-b).
\]
Bounded integer products use binary expansion. Added variables are existentially projected for
comparison. Sequential redundancy deletion retains the model after each proved deletion, and a
final relation query checks the reconstructed pair.

\section{Formal Guarantees and Complete Procedures}\label{app:formalguarantees}

This section establishes the certificate guarantees used by verification, slicing, reconstruction
and adoption. It then gives the complete procedures that instantiate those guarantees.

\subsection{Relation certificates}
\textbf{Proof of the relation certificate.} The two unsatisfiable difference queries establish both inclusions and equality
of the projected feasible regions. The preference query excludes any common decision pair ordered
strictly by one model and tied or oppositely ordered by the other. Their strict preference relations
are equal. A satisfying query supplies the separating point or decision pair directly.

A feasible model entails a compiled source obligation when its violation formula is unsatisfiable.
The returned result queries use the model, report and tolerance from Eqs.~\eqref{eq:solution}
and~\eqref{eq:better}. Extension and value binding are existential feasibility checks. The absence
of a better decision certifies optimality to the declared margin. An empty model receives its
separate infeasibility outcome.

\subsection{Constraint reduction and witness completion}
\textbf{Proof of sound reduction.} The full violation formula implies the reduced formula. The
reduced formula's unsatisfiability therefore implies the full formula's unsatisfiability. A complete
satisfying assignment lies in the candidate feasible region and outside the obligation. A partial
satisfying assignment is extended and checked against the full model before it supplies this witness.

\textbf{Proof of exact component slicing.} Omitted components share no variables with retained
components or the obligation. Base feasibility supplies a satisfying assignment in every omitted
component. Combining these assignments extends every satisfying retained assignment to the full
model. The reverse implication follows by restricting a full assignment to retained variables.

Neighborhood slices smaller than the component use Eq.~\eqref{eq:localrules} and full witness completion.
For a target with all variables observable, a directional difference decomposes into violations
of its individual rows and domains. Hidden auxiliaries retain the complete projected formula.
The evidence record stores the assignment, violated rows and relation used by the revision
operator in Section 4.4.

\subsection{Reconstruction and adoption}
\textbf{Proof of Proposition 1.} Let \(T\) be the selected target and
\(\mathcal N(T)\) its normalized representation. Preservation of the projected region and objective
graph gives \(F_{\mathcal N(T)}^X=F_T^X\) and \(G_{\mathcal N(T)}^X=G_T^X\). Compilation constructs
the normalized tuple, and runtime capture recovers it from the generated program. Hence
\[
F_{\widetilde B}^X=F_{\mathcal N(T)}^X=F_T^X,\qquad
G_{\widetilde B}^X=G_{\mathcal N(T)}^X=G_T^X.
\]
The exported mathematical model is \(\widetilde A=T\). Its region and objective graph agree with
those of the captured implementation. The final queries check both region inclusions and the
objective relation for these constructed artifacts.

Sequential redundancy deletion preserves equivalence by induction. At each step the remaining
rows entail the deleted row, and adding or removing that entailed row preserves the feasible region.
The objective is unchanged. Each new deletion uses the currently retained basis.

The adoption loop accepts a candidate only after its required relation checks succeed. After an
accepted construction, a failed proposal retains that certified pair and its evidence. Before such
acceptance, the return value retains the original verdicts. Source interpretation follows the compiled
requirements and retained provenance of the selected target.

\subsection{Construction and adoption procedures}\label{app:algorithmsteps}

Algorithm~\ref{alg:completion} states the complete verification, candidate generation and selection loop.
Algorithms~\ref{alg:constructiondetails} and~\ref{alg:acceptancedetails} expand its helper procedures.
The selection order gives the mathematical artifact priority when both inputs are readable.
Construction returns the mathematical model, regenerated program and computed report together.
The verification procedure applies the region, objective and report queries in Section~\ref{subsec:deciding2},
and checks the compiled requirements enabled by the source review configuration. The policy \(\rho\)
chooses direct, certified voting or protected selection without changing candidate verification.

\begin{algorithm}[htbp]
\caption{Verification and repair}\label{alg:completion}
\small
\begin{algorithmic}[1]
\Require Source \(s\), compiled requirements \(C\), artifacts \(Y_0=(A,\widehat P,R)\), budget \(K\), policy \(\rho\)
\Ensure Returned artifacts \(Y\) and verification evidence \(e_Y\)
\State \(e_0\gets\Call{Verify}{Y_0,C}\), \quad \(\mathcal B\gets\varnothing\)
\State \(T\gets\Call{SelectModel}{A,\widehat P}\)
\If{\(T\ne\bot\)}
    \State \(\mathcal B\gets\mathcal B\cup\{\Call{Complete}{T,C}\}\)
\EndIf
\For{\(k=1,\ldots,K\)}
    \State \(\widehat P'_k\gets\mathcal V(s,Y_0,e_0;k)\), \quad \(T'_k\gets\mathcal H(\widehat P'_k)\)
    \If{\(T'_k\ne\bot\)}
        \State \(\mathcal B\gets\mathcal B\cup\{\Call{Complete}{T'_k,C}\}\)
    \EndIf
\EndFor
\State \((Y,e)\gets\Call{SelectCertified}{\mathcal B,Y_0,\rho}\)
\If{\(Y=\bot\)}
    \State \Return \((Y_0,e_0)\)
\EndIf
\State \Return \((Y,e)\)
\end{algorithmic}
\end{algorithm}

\textbf{State and evidence.} Each proposal $Y$ and its evidence $e$ form one revision record. Direct
selection chooses the first fully certified record. Certified voting applies the value class quorum,
and protected selection first applies its declared preservation filter. If selection returns no
record, the original state and evidence remain active. The history records the selected target,
triggering verdict, proposed revision and verification after construction.

\begin{algorithm}[htbp]
\caption{Model selection and construction}\label{alg:constructiondetails}
\small
\begin{algorithmic}[1]
\Procedure{SelectModel}{$A,\widehat P$}
    \If{$A\neq\bot$}
        \State $T\gets A$
    \Else
        \State $T\gets\mathcal H(\widehat P)$
    \EndIf
    \State \Return $T$
\EndProcedure
\Statex
\Procedure{Complete}{$T,C$}
    \State $\widetilde A\gets T$, \quad
           $\widetilde P\gets\operatorname{Compile}(\mathcal N(T))$
    \State $\widetilde B\gets\mathcal H(\widetilde P)$, \quad
           $\widetilde R\gets\operatorname{Solve}(\widetilde P)$
    \State $Y\gets(\widetilde A,\widetilde P,\widetilde R)$
    \State $e\gets\Call{Verify}{Y,C}$
    \State \Return $(Y,e)$
\EndProcedure
\end{algorithmic}
\end{algorithm}

\begin{algorithm}[htbp]
\caption{Certified candidate acceptance}\label{alg:acceptancedetails}
\small
\begin{algorithmic}[1]
\Procedure{SelectCertified}{$\mathcal B,Y_0,\rho$}
    \State $E\gets\{(Y,e)\in\mathcal B:\text{ every relation required by }\rho\text{ is certified in }e\}$
    \If{$\rho$ is protected}
        \State $E\gets\{(Y,e)\in E:\Call{Preserves}{Y_0,Y,\rho}\}$
    \EndIf
    \If{$E=\varnothing$}
        \State \Return $(\bot,\varnothing)$
    \EndIf
    \If{$\rho$ uses certified voting}
        \State partition $E$ by exact equality of the certified value $v(Y)$
        \State $G\gets$ first largest value class in the fixed candidate order
        \If{$|G|<2$}
            \State \Return $(\bot,\varnothing)$
        \EndIf
        \State \Return the first $(Y,e)\in G$
    \EndIf
    \State \Return the first $(Y,e)\in E$ in fixed order
\EndProcedure
\end{algorithmic}
\end{algorithm}

\section{SMT Checking Coverage, Errors and Efficiency}\label{app:scale}
\subsection{Checking coverage and unresolved cases}
\begin{figure}[htbp]
\centering
\includegraphics[width=\textwidth]{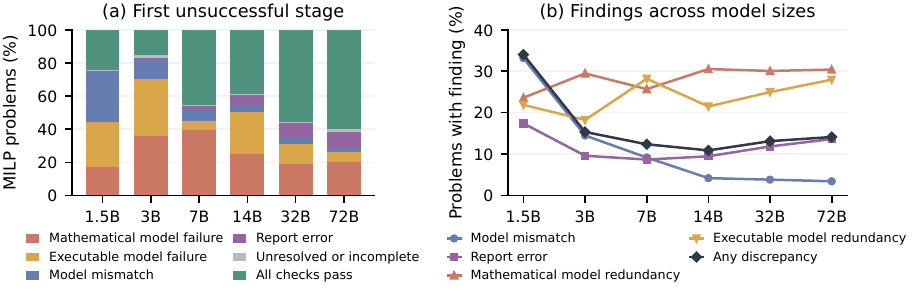}
\caption{Verification outcomes and detected findings over 3,034 problems per generator.}
\label{fig:stage1errors}
\end{figure}

Table~\ref{tab:ckptmain} separates artifact availability, model and report discrepancies, and the full
formal chain. Figure~\ref{fig:stage1errors} distinguishes the first unsuccessful stage from findings
that can overlap on the same output. Redundant constraints are diagnosed separately by entailment.
\begin{table}[H]
\caption{Original verification outcomes over 3,034 problems per generator, in percent.}
\label{tab:ckptmain}
\centering\footnotesize\setlength{\tabcolsep}{3pt}
\begin{tabular*}{\textwidth}{@{\extracolsep{\fill}}lrrrrrr@{}}
\toprule
Generator & \shortstack{Written model\\available} & \shortstack{Executable model\\available} & \shortstack{Model--code\\discrepancy} & \shortstack{Report\\discrepancy} & \shortstack{Any\\discrepancy} & \shortstack{All checks\\pass} \\
\midrule
1.5B & 85.6 & 64.3 & 33.2 & 17.4 & 34.0 & 24.2 \\
3B & 74.3 & 42.0 & 14.5 & 9.6 & 15.3 & 15.5 \\
7B & 72.6 & 75.8 & 9.2 & 8.7 & 12.4 & 45.2 \\
14B & 84.2 & 60.0 & 4.2 & 9.5 & 10.9 & 39.0 \\
32B & 89.2 & 77.3 & 3.8 & 11.8 & 13.1 & 55.8 \\
72B & 88.6 & 83.5 & 3.4 & 13.6 & 14.1 & 59.8 \\
\bottomrule
\end{tabular*}
\end{table}

\subsection{Input failures and model size}\label{app:firstlayeraudit}
Availability failures have one primary parsing or execution category, whereas semantic findings may
overlap. The following tables give the input failures behind the main availability rates and the
size of the models that reach verification. Individual traces and unresolved query identifiers
remain in the supplementary records.

\begin{table}[htbp]
\caption{Written-model input failures, with 3,034 assigned problems per generator.}
\label{tab:mathfailurecodes}
\centering\footnotesize\setlength{\tabcolsep}{2pt}
\begin{tabular*}{\textwidth}{@{\extracolsep{\fill}}lrrrrrr@{}}
\toprule
Written-model failure & 1.5B & 3B & 7B & 14B & 32B & 72B\\
\midrule
Malformed literal & 283 & 551 & 610 & 291 & 201 & 238\\
Undeclared symbol & 34 & 82 & 101 & 117 & 88 & 91\\
Missing expression & 21 & 54 & 31 & 10 & 29 & 12\\
Invalid comparator & 29 & 19 & 61 & 10 & 0 & 0\\
Missing model fields & 54 & 15 & 1 & 15 & 1 & 0\\
Missing right hand side & 0 & 33 & 7 & 0 & 0 & 0\\
Duplicate variable declaration & 6 & 4 & 11 & 4 & 3 & 0\\
Missing artifact & 4 & 4 & 4 & 4 & 4 & 4\\
Undeclared symbols & 3 & 1 & 0 & 20 & 0 & 0\\
Nonnumeric value & 0 & 7 & 0 & 0 & 0 & 0\\
Unrecognized objective fields & 3 & 1 & 2 & 1 & 0 & 0\\
Unrecognized constraint fields & 0 & 5 & 1 & 0 & 0 & 0\\
Ambiguous term fields & 0 & 0 & 0 & 4 & 1 & 0\\
Missing left hand side & 0 & 4 & 0 & 0 & 0 & 0\\
Invalid operator arity & 0 & 0 & 0 & 3 & 0 & 0\\
Duplicate or invalid key & 0 & 0 & 1 & 0 & 0 & 0\\
Invalid variable domain & 1 & 0 & 0 & 0 & 0 & 0\\
Nonliteral construct & 0 & 1 & 0 & 0 & 0 & 0\\
Unbound or unsupported expression & 0 & 0 & 1 & 0 & 0 & 0\\
Unrecognized variable fields & 0 & 0 & 0 & 1 & 0 & 0\\
\midrule all unavailable & 438 & 781 & 831 & 480 & 327 & 345
\\
\bottomrule
\end{tabular*}
\end{table}

\begin{table}[htbp]
\caption{Executable-model input failures, with 3,034 assigned problems per generator.}
\label{tab:codefailurecodes}
\centering\footnotesize\setlength{\tabcolsep}{2pt}
\begin{tabular*}{\textwidth}{@{\extracolsep{\fill}}lrrrrrr@{}}
\toprule
Executable-model failure & 1.5B & 3B & 7B & 14B & 32B & 72B\\
\midrule
Python runtime error & 811 & 1,224 & 211 & 865 & 442 & 170\\
No program in response & 239 & 444 & 434 & 229 & 175 & 225\\
Python syntax error & 28 & 86 & 77 & 110 & 67 & 93\\
Truncated before solve & 0 & 1 & 8 & 7 & 2 & 8\\
Missing artifact & 4 & 4 & 4 & 4 & 4 & 4\\
\midrule all unavailable & 1,082 & 1,759 & 734 & 1,215 & 690 & 500
\\
\bottomrule
\end{tabular*}
\end{table}

\begin{table}[htbp]
\caption{Generated artifact size by generator, with 3,034 assigned problems per generator.}
\label{tab:stage1sizes}
\centering\footnotesize\setlength{\tabcolsep}{2pt}
\begin{tabular*}{\textwidth}{@{\extracolsep{\fill}}llrrrrrr@{}}
\toprule
Generator & Artifact & Assigned & Available & \shortstack{Median\\variables} & \shortstack{Median\\constraints} & \shortstack{Maximum\\variables} & \shortstack{Maximum\\constraints}\\
\midrule
1.5B & Written model & 3,034 & 2,596 & 2 & 1 & 81 & 100\\
1.5B & Executable model & 3,034 & 1,952 & 3 & 2 & 450 & 96\\
3B & Written model & 3,034 & 2,253 & 3 & 3 & 79 & 130\\
3B & Executable model & 3,034 & 1,275 & 2 & 3 & 650 & 1,913\\
7B & Written model & 3,034 & 2,203 & 3 & 3 & 98 & 119\\
7B & Executable model & 3,034 & 2,300 & 3 & 3 & 336 & 386\\
14B & Written model & 3,034 & 2,554 & 3 & 4 & 156 & 134\\
14B & Executable model & 3,034 & 1,819 & 4 & 4 & 15,150 & 23,396\\
32B & Written model & 3,034 & 2,707 & 3 & 4 & 201 & 116\\
32B & Executable model & 3,034 & 2,344 & 3 & 3 & 1,891 & 3,821\\
72B & Written model & 3,034 & 2,689 & 3 & 3 & 128 & 171\\
72B & Executable model & 3,034 & 2,534 & 3 & 3 & 5,150 & 10,070
\\
\bottomrule
\end{tabular*}
\end{table}

Availability is counted separately for written and executable models because either artifact may be
missing while the other remains usable. Medians describe the typical captured artifact, whereas the
maxima expose the larger models that determine the checking limits. Table~\ref{tab:formalstructures}
then separates affine models, other admitted structures and unavailable inputs on the same assignments.

\begin{table}[htbp]
\caption{Admitted structures and unavailable inputs over 3,034 assigned problems per generator.}
\label{tab:formalstructures}
\centering\footnotesize\setlength{\tabcolsep}{2pt}
\begin{tabular*}{\textwidth}{@{\extracolsep{\fill}}lrrrrrrr@{}}
\toprule
& Assigned & \multicolumn{3}{c}{Written model} & \multicolumn{3}{c}{Executable model}\\
\cmidrule(lr){3-5}\cmidrule(lr){6-8}
Model size & & Affine & \shortstack{Other admitted\\structure} & Unavailable & Affine & \shortstack{Other admitted\\structure} & Unavailable\\
\midrule
1.5B & 3,034 & 2,595 & 1 & 438 & 1,949 & 3 & 1,082\\
3B & 3,034 & 2,249 & 4 & 781 & 1,265 & 10 & 1,759\\
7B & 3,034 & 2,202 & 1 & 831 & 2,297 & 3 & 734\\
14B & 3,034 & 2,551 & 3 & 480 & 1,811 & 8 & 1,215\\
32B & 3,034 & 2,695 & 12 & 327 & 2,334 & 10 & 690\\
72B & 3,034 & 2,689 & 0 & 345 & 2,527 & 7 & 500
\\
\bottomrule
\end{tabular*}
\end{table}

\subsection{Errors introduced during code generation}\label{app:errorattribution}

\begin{table}[htbp]
\caption{Reusable finding taxonomy. Each label is tied to checkable evidence rather than inferred from a final objective value.}
\label{tab:findingtaxonomy}
\centering\footnotesize\setlength{\tabcolsep}{2pt}
\begin{tabular*}{\textwidth}{@{\extracolsep{\fill}}p{0.22\textwidth}p{0.72\textwidth}@{}}
\toprule
Finding & Required evidence \\
\midrule
Additional or missing restriction & A separating assignment is feasible on one side and excluded by a named row on the other. \\
Domain or bound difference & A separating assignment violates a declaration while satisfying the target rows. \\
Objective value or preference difference & A shared feasible point has different values, or two shared points receive different orderings. \\
Redundant constraint & The remaining rows and fixed domains entail the row under sequential deletion. \\
Report error & The reported point has no feasible completion, cannot bind the reported value, or admits a better feasible point. \\
Transcription error & Aligned rows or declarations differ in coefficient, index, bound or direction, with source provenance. \\
\bottomrule
\end{tabular*}
\end{table}

The executable model introduces a distinct and larger failure surface even though both artifacts
originate from the same model response. Across 3,034 problems and six generators (18,204 outputs), code is unavailable in 5,980
cases (32.9\%), compared with 3,202 unavailable written models (17.6\%). On the 951 pairs for which
both sequential constraint bases are completely decided, the code adds at least one restriction not
entailed by the written model in 546 cases and omits at least one written restriction in 415 cases;
these categories overlap. Across the six generator sizes, redundancy appears only in code in 239
assignments and only in the written model in 188.  Table~\ref{tab:findingtaxonomy} defines the evidence required for each
reusable finding. This asymmetry motivates checking the generated implementation itself rather than
treating successful execution or a matching objective value as evidence that transcription preserved
the formulation.

\subsection{Checking time on the controlled CPU panel}\label{app:stage1cost}
The controlled experiment selects the available executable model with the most constraints for
each problem, breaking ties by variables and then generator size. It runs on two AMD EPYC 9555
processors with 128 physical cores and 755 GiB of memory, using eight concurrent problems.
Each execution mode has a ten second total deadline covering input loading, preprocessing, solver
processes and evidence writing. Historical evidence caches are disabled. Median checked time includes
timeouts. Running both execution modes takes 541.9 seconds. Problems without a captured executable
model remain unresolved in the 3,034-problem denominator.

\begin{figure}[htbp]
\centering
\includegraphics[width=0.98\textwidth]{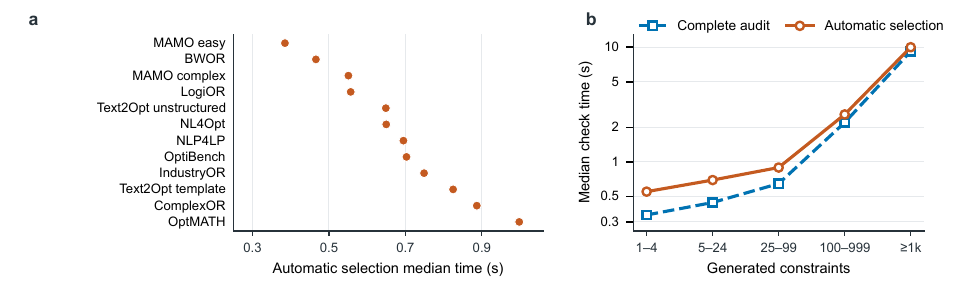}
\caption{ORSMT checking time across (a) benchmarks with automatic backend selection and (b) constraint sizes for the complete audit and automatic selection. The panel contains all 3,034 problems. Unavailable models are unresolved, and timeouts enter the medians.}
\label{fig:smtscaling}
\end{figure}

The complete audit decides 2,860 problems with a 0.397 second median. Automatic backend selection
decides 2,869 problems with a 0.647 second median. Both settings run the same ORSMT relation checks
and agree on every problem that both decide. Table~\ref{tab:coldfamilies} gives automatic selection
results by benchmark. Figure~\ref{fig:smtscaling} gives the corresponding constraint size comparison.

\begin{table}[htbp]
\caption{ORSMT checking by benchmark on the 3,034-problem CPU panel.}
\label{tab:coldfamilies}
\centering\footnotesize\setlength{\tabcolsep}{2pt}
\begin{tabular*}{\textwidth}{@{\extracolsep{\fill}}lrrrrr@{}}
\toprule
Benchmark & Problems & Decided & Unresolved & \shortstack{Problems\\decided (\%)} & \shortstack{Median checked\\time (s)}\\
\midrule
MAMO easy & 652 & 652 & 0 & 100.0 & 0.385\\
Text2Opt template & 630 & 502 & 128 & 79.7 & 0.825\\
OptiBench & 422 & 422 & 0 & 100.0 & 0.703\\
Text2Opt unstructured & 248 & 228 & 20 & 91.9 & 0.649\\
NLP4LP & 242 & 242 & 0 & 100.0 & 0.695\\
NL4Opt & 230 & 230 & 0 & 100.0 & 0.650\\
MAMO complex & 211 & 211 & 0 & 100.0 & 0.551\\
OptMATH & 132 & 122 & 10 & 92.4 & 0.998\\
IndustryOR & 98 & 95 & 3 & 96.9 & 0.749\\
LogiOR & 76 & 73 & 3 & 96.1 & 0.557\\
BWOR & 75 & 74 & 1 & 98.7 & 0.466\\
ComplexOR & 18 & 18 & 0 & 100.0 & 0.887\\
\midrule
All benchmarks & 3,034 & 2,869 & 165 & 94.6 & 0.647
\\
\bottomrule
\end{tabular*}
\end{table}

\begin{table}[htbp]
\caption{ORSMT checking by variable domain over the complete 3,034-problem panel.}
\label{tab:coldtypes}
\centering\footnotesize\setlength{\tabcolsep}{2pt}
\begin{tabular*}{\textwidth}{@{\extracolsep{\fill}}lrrrr@{}}
\toprule
Variable domain & Problems & Decided & \shortstack{Problems\\decided (\%)} & \shortstack{Median checked\\time (s)}\\
\midrule
No captured model & 134 & 0 & 0.0 & --\\
Binary & 119 & 117 & 98.3 & 0.689\\
Continuous & 406 & 404 & 99.5 & 0.700\\
Integer & 1,951 & 1,942 & 99.5 & 0.598\\
Mixed integer & 424 & 406 & 95.8 & 0.847\\
\midrule
All problems & 3,034 & 2,869 & 94.6 & 0.647
\\
\bottomrule
\end{tabular*}
\end{table}

\subsection{Canonical requirements and redundancy}\label{app:canonicalpanel}
Canonical families expose exact requirement queries and row entailment at different model
structures. The timing table uses one model as its unit. The redundancy table uses one tested row
and records unresolved queries separately. Sequential removal of entailed rows produces a smaller
equivalent basis for subsequent comparisons and reconstruction.

\begin{table}[htbp]
\caption{Exact requirement queries over 136,100 obligations in 1,806 canonical OR models. Times are seconds per model.}
\label{tab:familyverification}
\centering\footnotesize\setlength{\tabcolsep}{2pt}
\begin{tabular*}{\textwidth}{@{\extracolsep{\fill}}lrrrr@{}}
\toprule
Family & Models & Obligations & Median rows & Median (s)\\
\midrule
BinPacking & 620 & 31,000 & 399 & 0.257\\
Inventory & 160 & 12,000 & 75 & 0.013\\
JobShop & 160 & 60,320 & 280 & 0.063\\
Knapsack & 260 & 260 & 1 & 0.005\\
NetworkFlow & 50 & 3,250 & 63 & 0.014\\
Pollution & 80 & 600 & 8 & 0.015\\
Portfolio & 160 & 4,640 & 29 & 0.007\\
Production & 80 & 520 & 6 & 0.004\\
TSP & 156 & 22,470 & 123 & 0.052\\
Transportation & 80 & 1,040 & 13 & 0.009\\
\bottomrule
\end{tabular*}
\end{table}

\begin{table}[htbp]
\caption{Exact row entailment over 340,289 rows in 1,705 canonical OR models.}
\label{tab:canonicalredundancy}
\centering\footnotesize\setlength{\tabcolsep}{2pt}
\begin{tabular*}{\textwidth}{@{\extracolsep{\fill}}lrrrrr@{}}
\toprule
Family & Checked & Rows & Redundant & Unresolved & Models with redundancy\\
\midrule
BinPacking & 575 & 240,761 & 2,336 & 2,612 & 575 \\
Inventory & 160 & 12,000 & 3,355 & 0 & 160 \\
JobShop & 155 & 55,320 & 0 & 2,205 & 0 \\
Knapsack & 260 & 260 & 0 & 0 & 0 \\
NetworkFlow & 50 & 3,250 & 1,803 & 0 & 50 \\
Pollution & 34 & 242 & 0 & 0 & 0 \\
Portfolio & 155 & 4,426 & 1,912 & 0 & 155 \\
Production & 80 & 520 & 45 & 0 & 38 \\
TSP & 156 & 22,470 & 3,586 & 15 & 156 \\
Transportation & 80 & 1,040 & 0 & 0 & 0
\\
\midrule
Total & 1,705 & 340,289 & 13,037 & 4,832 & 1,134\\
\bottomrule
\end{tabular*}
\end{table}

\subsection{Constraint selection and larger instances}\label{app:executablescaling}
The FrontierOR scaling panel~\citep{kong2026frontieror} compares complete and restricted queries on
the same instantiated models. Table~\ref{tab:frontierslicingbysize} reports the paired timing
population.

\begin{table}[htbp]
\caption{Timing comparison on 5,869 FrontierOR claims below 100,000 constraints.}
\label{tab:frontierslicingbysize}
\centering\footnotesize\setlength{\tabcolsep}{2pt}
\begin{tabular*}{\textwidth}{@{\extracolsep{\fill}}lrrrrrr@{}}
\toprule
Constraints & Claims & \shortstack{Complete\\queries} & \shortstack{Decided by\\both} & \shortstack{Complete\\median (s)} & \shortstack{Restricted\\median (s)} & \shortstack{Paired\\speedup}\\
\midrule
Under 1,000 & 3,655 & 3,655 & 3,589 & 0.79 & 0.29 & 1.9\\
1,000 to 10,000 & 1,315 & 1,315 & 1,237 & 11.79 & 0.30 & 29.3\\
10,000 to 100,000 & 899 & 249 & 229 & 44.71 & 1.30 & 27.8\\
\bottomrule
\end{tabular*}
\end{table}

A reduced unsatisfiable query proves entailment. Candidate refutations are completed with a native
solver and checked against every model row using tolerance
\(10^{-6}\max(1,|b|)+10^{-6}\) for right hand side \(b\). Local entailment uses exact arithmetic.
The construction cap limits full queries in the largest interval. All 5,055 claims decided by both
query types agree.

The extension beyond 20,000 rows deduplicates query attempts by family, instance, row and edit
kind, retaining the latest record. It contains 1,383 unique claims on 78 models and excludes base
status records. Each full model is constructed once and reused with push and pop. Capture, encoding,
base feasibility and query time have separate records.

\begin{table}[htbp]
\caption{Complete query coverage for 1,383 unique FrontierOR claims beyond 20,000 constraints.}
\label{tab:abovecap}
\centering\footnotesize\setlength{\tabcolsep}{2pt}
\begin{tabular*}{\textwidth}{@{\extracolsep{\fill}}lrrrrr@{}}
\toprule
Model rows & Models & Claims & Full query decisions & Checks decided (\%) & Matching decisions\\
\midrule
20,000 to 50,000 & 9 & 221 & 115 & 52.0 & 112/112 \\
50,000 to 100,000 & 23 & 568 & 468 & 82.4 & 439/439 \\
100,000 to 1 million & 41 & 574 & 75 & 13.1 & 49/49 \\
At least 1 million & 5 & 20 & 0 & 0.0 & 0/0
\\
\bottomrule
\end{tabular*}
\end{table}

Full queries decide 658 claims, with matching verdicts on all 600 claims decided by both query types.
Exact base feasibility is available for 584 claims, of which 537 full queries decide. Among
799 claims with unresolved base feasibility, 121 full queries decide. The largest recorded models
reach 5,768,230 rows and 3,302,127 variables. One 372,000-row model takes 283 seconds to encode and
108 seconds to establish an exact feasible base. A weakened row query takes 5,656 seconds in the
complete solver and 2.60 seconds locally. These are measurements of that individual query.

An independent solver comparison constructs obligations separately for Z3 and cvc5. It records
377,151 of 377,280 assigned query pairs, with agreement on all 110,366 pairs decided by both solvers.
At least one solver remains unresolved on 266,785 recorded pairs, and 129 pairs lack a completed
record. This experiment checks independent solver encodings and verdicts.

\section{Feedback Quality, Repair Effects and Transfer}\label{app:repair}
This section separates three questions. First, what does each verifier check, and what evidence does
it return? Second, how much do direct reconstruction and language-model revision repair under the
same data and budget? Third, do the repaired pipelines improve correct answers, model--code matches,
source requirements and host integrations without hiding residual constraint differences?

\subsection{Feedback coverage and error localization}
ReLoop evaluates execution behavior and selected perturbations, whereas Opt-Verifier produces
structure and solution judgments~\citep{lian2026reloop,liu2026optverifier}. Table~\ref{tab:verifiercoveragefull}
reports their available decisions and flags on the exact model--code cohorts. Because these checks
target different relations, the table keeps model--code refutations and matches in separate columns.
\begin{table}[htbp]
\caption{Feedback coverage over 3,034 problems and six generators (18,204 outputs), with refuted and equivalent model--code pairs separated.}
\label{tab:verifiercoveragefull}
\centering\footnotesize\setlength{\tabcolsep}{1.5pt}
\begin{tabular*}{\textwidth}{@{\extracolsep{\fill}}lrrrrr@{}}
\toprule
Verifier & \shortstack{Outputs with\\a decision\\/ 18,204} & \shortstack{Refuted pairs\\flagged\\/ 1,905} & \shortstack{Refuted pairs\\not flagged\\/ 1,905} & \shortstack{No decision on\\refuted pairs\\/ 1,905} & \shortstack{Equivalent pairs\\flagged\\/ 7,867} \\
\midrule
ReLoop execution & 16,331 & 181 & 1,716 & 8 & 29 \\
ReLoop perturbation & 803 & 324 & 51 & 1,530 & 369 \\
Opt-Verifier structure & 15,378 & 1,656 & 249 & 0 & 4,576 \\
Opt-Verifier solution & 10,202 & 1,633 & 108 & 164 & 7,001 \\
\bottomrule
\end{tabular*}
\end{table}

FLARE verifies predefined formulation-family relations~\citep{robbins2026flare}. Its unit is therefore
different from the instance-level mathematical model, executable model and report studied here. We
use FLARE to distinguish the scope of the verification task; the paired feedback and repair study uses
ReLoop and Opt-Verifier on the shared 1,905-output cohort.

ORSMT returns either a certificate or a separating assignment tied to rows, domains or objective
preferences. Two checked cases expose errors in comparison method feedback. ReLoop accepts the
IndustryOR 71 result although the mathematical and implemented feasible regions differ in both
directions. Opt-Verifier rejects the certified BWOR 18 value 34 using an alternative whose listed
costs sum to 37. Appendix~\ref{app:workedexamples} gives the complete models, arithmetic and repair
evidence.

\subsection{Direct reconstruction and completion routes}
Table~\ref{tab:relooppasstypes} decomposes the ReLoop execution passes by overlapping semantic
findings.
\begin{table}[htbp]
\caption{ReLoop passes by discrepancy category on 1,905 refuted outputs. Categories overlap.}
\label{tab:relooppasstypes}
\centering\footnotesize\setlength{\tabcolsep}{3pt}
\begin{tabular*}{\textwidth}{@{\extracolsep{\fill}}lrrr@{}}
\toprule
Discrepancy category & Affected outputs & ReLoop decisions & ReLoop passes \\
\midrule
Additional feasible decisions & 1,114 & 1,109 & 1,090 \\
Excluded mathematical decisions & 1,232 & 1,228 & 1,056 \\
Objective preference disagreement & 203 & 202 & 200 \\
Variable declaration differences & 922 & 915 & 836 \\
\bottomrule
\end{tabular*}
\end{table}

Tables~\ref{tab:completionroutes} and~\ref{tab:completionstates} distinguish the available source
artifact and the original failure. A mathematical target reconstructs the program, and a captured
program reconstructs the mathematics. This mechanical stage uses no language model calls.
\begin{table}[htbp]
\caption{Mechanical completion routes over 3,034 problems and six generators (18,204 outputs).}
\label{tab:completionroutes}
\centering\footnotesize\setlength{\tabcolsep}{3pt}
\begin{tabular*}{\textwidth}{@{\extracolsep{\fill}}lrrr@{}}
\toprule
Starting artifact & Assigned outputs & \shortstack{Model--code\\matches} & \shortstack{Including empty\\feasible regions} \\
\midrule
Already certified pair & 7,867 & 7,867 & 0 \\
Readable written model & 6,590 & 6,590 & 1,396 \\
Captured executable model & 1,267 & 1,267 & 321 \\
Neither artifact usable & 2,480 & 0 & 0 \\
\bottomrule
\end{tabular*}
\end{table}

\begin{table}[htbp]
\caption{Mechanical completion outcomes on 10,337 initially uncertified outputs.}
\label{tab:completionstates}
\centering\footnotesize\setlength{\tabcolsep}{3pt}
\begin{tabular*}{\textwidth}{@{\extracolsep{\fill}}lrrr@{}}
\toprule
Initial failure & Assigned outputs & Usable counterpart & \shortstack{Model--code\\matches} \\
\midrule
Original executable-model error & 3,319 & 3,319 & 3,319 \\
Missing model or implementation & 2,826 & 839 & 839 \\
Proved model inconsistency & 1,905 & 1,905 & 1,905 \\
Unsupported written representation & 1,092 & 599 & 599 \\
Unresolved relation queries & 904 & 904 & 904 \\
No shared decision variables & 174 & 174 & 174 \\
Only some variables matched & 72 & 72 & 72 \\
Unsupported runtime capture & 25 & 25 & 25 \\
Program does not invoke optimization & 16 & 16 & 16 \\
Comparison error & 4 & 4 & 4 \\
\bottomrule
\end{tabular*}
\end{table}

Table~\ref{tab:completionmatrix} in the main paper retains every benchmark and generator combination.
Protected selection retains finite original results, reconstructs paired artifacts around those
values and applies source revision to outputs without a finite result.
\begin{figure}[htbp]
\centering
\includegraphics[width=\textwidth]{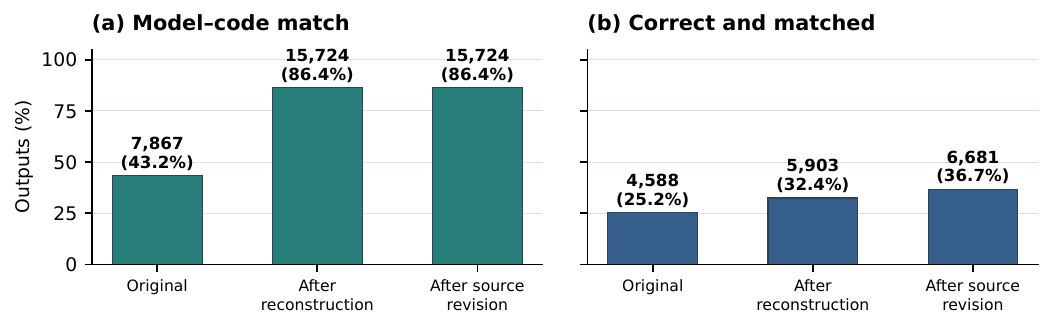}
\caption{Model--code matches and outputs that are both correct and matched before protected selection over 3,034 problems and six generators (18,204 outputs).}
\label{fig:repairoutcomes}
\end{figure}

\subsection{LLM feedback, ablations and paired repair}\label{app:completiondetails}
Protected selection retains every finite original result. It reconstructs a paired artifact around
that result and applies value changing revision when the original has no finite value. The
mathematical artifact is the default reconstruction target, and the captured program supplies the
target when the mathematical artifact is unavailable. Every route records the selected target.

The fixed 14B reviser receives the problem statement and evidence for its repair condition. All
single candidate revision calls use temperature zero with a 32,768 token context. Program, model, requirement
and judgment allowances are 3,072, 4,096, 8,192 and 800 tokens, respectively. Native optimization
receives 30 seconds, rational reconstruction five seconds, and final checks ten seconds per query
with 120 seconds per relation. The five candidate policy groups certified finite values by exact
equality, select the largest group and break ties by the fixed candidate order. Reference values enter
scoring after selection.

\subsubsection{Prompts used for language model revision}
The matched experiment compares direct revision, which receives only a generic mismatch notice,
with counterexample-guided revision, which additionally receives the failed relation and separating
assignment. Both conditions use the same system instruction, problem statement, written model,
previous program, reviser and decoding budget. Bracketed fields below are filled from the stored
record for that output. Neither condition receives a benchmark answer or reference objective value.

\begin{evidencecase}{Prompt for direct revision}
\footnotesize\raggedright
\textbf{System:} You are an operations-research modeling expert.

\smallskip
\textbf{User:} You wrote a mathematical model and a gurobipy program for the problem below, and
the program does not implement the model you wrote. Return corrected gurobipy code only, with no
prose or Markdown fences. Change only what is wrong.

\smallskip
\textbf{Problem:} [problem statement]

\textbf{Your stated model:} [written model in JSON]

\textbf{Your program:} [previous program]

\textbf{Repair instruction:} A check found that the program and the model do not define the same
problem.
\end{evidencecase}

\begin{evidencecase}{Prompt for counterexample-guided revision}
\footnotesize\raggedright
\textbf{System:} You are an operations-research modeling expert.

\smallskip
\textbf{User:} You wrote a mathematical model and a gurobipy program for the problem below. The
written model is the implementation target. Correct only the program.

\smallskip
\textbf{Problem:} [problem statement]

\textbf{Written model (JSON):} [complete written model]

\textbf{Previous program:} [complete previous program]

\textbf{Verification evidence:} A checker found a concrete disagreement. [failed relations]
[written-model constraints involved] [runtime-program constraints involved] [variables involved]
[an assignment accepted by the program but forbidden by the model] [an assignment allowed by the
model but rejected by the program] [two assignments ranked differently] [domain or bound differences]

\textbf{Repair instruction:} Use the counterexamples to correct the program. Return complete
corrected code only.
\end{evidencecase}

\begin{evidencecase}{Comparison-feedback prompt variants}
\footnotesize\raggedright
Both variants use the system instruction, problem statement, written model, previous program and
code-only output instruction shown above. Only the following verifier-feedback suffix changes.

\smallskip
\textbf{ReLoop suffix:} A verifier executed your program and reported: [one line per finding in the
form \texttt{check (severity): message}].

\smallskip
\textbf{Opt-Verifier suffix:} A verifier read your model and your reported decision and said:
\texttt{structure side:} [structure comment]; \texttt{solution side:} [solution comment].
\end{evidencecase}

Thus the comparison requests differ in verifier feedback rather than in the task or output format.
Every returned program is recaptured and checked before it can be selected.

The implementation serializes at most 6,000 characters of the written model and 8,000 characters of
the previous program. Each separating assignment is limited to 1,500 characters; the two assignments
in an objective-preference witness are limited to 700 characters each. ReLoop and Opt-Verifier
comments are limited to 4,000 characters. These limits and the code extractor are identical across
the matched repair conditions. The extractor accepts plain code or removes one Markdown code fence,
after which the complete returned program is executed, recaptured and checked.

Plan and row review start from the same source revision. Row review revises flagged candidates,
adopts finite certified proposals that pass recheck, and reuses the plan condition's infeasible
and unbounded branches. The compiled fallback retains previously certified pairs.
Table~\ref{tab:newrepairquality} reports the voting ablation. Table~\ref{tab:pairednew} compares its
paired outcomes with the two baselines. The interval resamples all outputs of each of the 1,531
source problems together for 10,000 draws with seed 20260913.

The one-round comparison gives every method the same refuted outputs and 14B reviser. The ORSMT arm
applies one source revision request followed by exact reconstruction. The comparison measures
model--code matches between the returned written and executable artifacts;
Table~\ref{tab:alignmentcost} reports the outcome and generation cost.
\begin{table}[htbp]
\caption{One-round model--code repair on 1,905 refuted outputs with Qwen2.5-14B.}
\label{tab:alignmentcost}
\centering\footnotesize\setlength{\tabcolsep}{0.9pt}
\begin{tabular*}{\textwidth}{@{\extracolsep{\fill}}lrrrrrr@{}}
\toprule
Repair method & \shortstack{Model--code\\matches} & \shortstack{Refuted pairs\\repaired (\%)} & \multicolumn{2}{c}{Per output} & \multicolumn{2}{c}{Per matched output}\\
\cmidrule(lr){4-5}\cmidrule(lr){6-7}
 & & & Calls & Tokens & Calls & Tokens\\
\midrule
ORSMT source revision and reconstruction & \textbf{1,607} & \textbf{84.4} & 1.00 & 3,596 & \textbf{1.19} & \textbf{4,263}\\
ReLoop & 109 & 5.7 & 1.00 & 2,830 & 17.5 & 49,460\\
Opt-Verifier & 68 & 3.6 & 5.80 & 7,949 & 162.5 & 222,690\\
\bottomrule
\end{tabular*}
\end{table}

\begin{table}[htbp]
\caption{Clustered bootstrap comparison on 1,905 refutations for five-candidate certified voting.}
\label{tab:pairednew}
\centering\footnotesize\setlength{\tabcolsep}{3pt}
\begin{tabular*}{\textwidth}{@{\extracolsep{\fill}}lrrrr@{}}
\toprule
Comparator & \shortstack{Voting\\only} & \shortstack{Comparator\\only} & \shortstack{Net gain\\percentage points} & \shortstack{95\% confidence interval\\percentage points} \\
\midrule
ReLoop & 175 & 59 & 6.09 & 4.54 to 7.66 \\
Opt-Verifier & 140 & 113 & 1.42 & $-0.26$ to 3.05 \\
\bottomrule
\end{tabular*}
\end{table}

The interval is above zero against ReLoop and overlaps zero against Opt-Verifier. Together with the
mechanical completion results, these comparisons separate repair without an LLM from revision driven
by verifier feedback and from multi-candidate selection.

\subsection{LLM repair cost and recorded latency}\label{app:auditcost}
Mechanical artifact reconstruction uses no language model calls. Source revision and candidate
selection add generation work after that step. Table~\ref{tab:generationcost} reports completed
records and recorded durations. Each duration measures one operation, and stages are reported
separately.

\begin{table}[htbp]
\caption{Recorded generation and repair work over 1,905 assigned outputs per condition.}
\label{tab:generationcost}
\centering\footnotesize\setlength{\tabcolsep}{3pt}
\begin{tabular*}{\textwidth}{@{\extracolsep{\fill}}lrrrr@{}}
\toprule
Repair setting & \shortstack{Candidate\\records} & \shortstack{Timed\\records} & \shortstack{Summed time\\(s)} & \shortstack{Mean time\\(s)} \\
\midrule
Source revision, 14B, one candidate & 1,905 & 1,905 & 107,068.3 & 56.20 \\
Source revision, 14B, five candidates & 9,525 & 9,525 & 425,946.3 & 44.72 \\
Source revision, 32B, one candidate & 1,905 & 1,905 & 260,619.2 & 136.81 \\
ReLoop revision, 14B & 1,905 & 1,903 & 35,219.4 & 18.51 \\
Opt-Verifier revision, 14B & 1,905 & 1,905 & 39,015.5 & 20.48 \\
\bottomrule
\end{tabular*}
\end{table}

Five-candidate repair contains one deterministic candidate and four samples per refutation.
Its summed generation duration corresponds to 223.59 recorded seconds per intervention, before
composition and verification. The 14B single-candidate mean is 56.20 seconds, and the 32B mean is
136.81 seconds on the same assigned cohort. Source-review work is tracked separately from candidate
generation, and conditional branches reuse existing candidates.

\subsubsection{Request and review accounting}\label{app:costdetails}

\textbf{Latency measurement.} Each candidate record stores one generated program and its duration.
The deterministic 14B and 32B conditions each contain one candidate per assigned output. The two
additional 14B files each contain two samples per output, giving four samples beyond the
deterministic candidate. Baseline records cover their respective revision operations. The duration
starts and ends at the boundaries of the recorded operation. Summing it across records measures
recorded work, while concurrent scheduling determines experiment wall time.

\begin{table}[htbp]
\caption{Generation duration over 15,238 timed candidate and baseline records. Times are seconds per operation.}
\label{tab:costquantiles}
\centering\footnotesize\setlength{\tabcolsep}{3pt}
\begin{tabular*}{\textwidth}{@{\extracolsep{\fill}}lrrr@{}}
\toprule
Operation & Timed records & Mean (s) & Median (s) \\
\midrule
14B deterministic candidate & 1,905 & 56.2 & 26.6 \\
14B samples 1 and 2 & 3,810 & 41.9 & 21.2 \\
14B samples 3 and 4 & 3,810 & 41.8 & 21.0 \\
32B deterministic candidate & 1,905 & 136.8 & 52.7 \\
ReLoop 14B revision & 1,903 & 18.5 & 9.0 \\
Opt-Verifier 14B revision & 1,905 & 20.5 & 10.1 \\
\bottomrule
\end{tabular*}
\end{table}

The means exceed the medians in every condition, and the deterministic 32B candidate has the largest
central duration on the same assigned outputs.

\textbf{Candidate accounting.} The five candidate policy evaluates one deterministic program and
four sampled programs for each refuted output. It reconstructs each usable program, checks the
resulting pair, solves the candidate model and groups eligible candidates by their certified values.
The selected candidate then enters the same adoption and fallback rule as the other conditions.
The deterministic candidate and the four samples account for all 9,525 candidate records. The
recorded time per intervention is the summed duration divided by the 1,905 assigned outputs.

\textbf{Review accounting.} Source review first checks the selected mathematical interpretation.
A proposed revision adds a generation operation, and a recheck evaluates the revised requirements.
Both approaches share 1,452 initial checks. Plan review adds 188 revision operations and 71
rechecks, while row review adds 783 revision operations and 716 rechecks. Row review triggers more
revisions because it tests individual source requirements that plan review can leave unexamined.
The returned revision becomes eligible for adoption after the formal checks and the source recheck.

\subsection{Full before-and-after repair results}\label{app:nativefull}
We now move from individual repair routes to complete before-and-after pipelines.
The main comparison assigns each of six generators the same 3,034 problems, for
$3{,}034\times6=18{,}204$ outputs.
Each method uses its own repair and adoption pipeline; ORSMT uses protected selection.
ORSMT produces 7,309 correct answers, 16,053 model--code matches and 7,308 outputs satisfying both. It
retains all 4,931 originally correct answers and adds 2,378. Correct answer counts and counts satisfying both use 18,192
numeric reference outputs, while model--code matches use all 18,204 assigned outputs from the
3,034 problems. Tables~\ref{tab:nativereferencebygenerator},
\ref{tab:nativecostdetail} and~\ref{tab:nativepairedcounts} give generator results, costs and paired discordant counts.
\begin{table}[htbp]
\caption{Correct answers by generator under protected selection over 18,192 outputs with numeric references.}
\label{tab:nativereferencebygenerator}
\centering\footnotesize\setlength{\tabcolsep}{3pt}
\begin{tabular*}{\textwidth}{@{\extracolsep{\fill}}lrrrrrrr@{}}
\toprule
Pipeline & 1.5B & 3B & 7B & 14B & 32B & 72B & Total\\
\midrule
Original & 213 & 278 & 801 & 860 & 1,326 & 1,453 & 4,931\\
ReLoop & 440 & 682 & 901 & 912 & 1,348 & 1,477 & 5,760\\
Opt-Verifier & \textbf{1,199} & 972 & 1,103 & 966 & 1,321 & 1,509 & 7,070\\
\textbf{ORSMT} & 684 & \textbf{1,037} & \textbf{1,163} & \textbf{1,378} & \textbf{1,506} & \textbf{1,541} & \textbf{7,309}\\
\bottomrule
\end{tabular*}
\end{table}

\begin{table}[htbp]
\caption{New and retained correct answers with language model work under protected selection over 3,034 problems and six generators (18,204 outputs).}
\label{tab:nativecostdetail}
\centering\footnotesize\setlength{\tabcolsep}{2.8pt}
\begin{tabular*}{\textwidth}{@{\extracolsep{\fill}}lrrrr@{}}
\toprule
Pipeline & \shortstack{New correct\\answers} & \shortstack{Original correct\\answers kept} & \shortstack{Language model\\calls} & \shortstack{Estimated tokens\\(millions)}\\
\midrule
Original & 0 & 4,931/4,931 & 0 & 0.0\\
ReLoop & 829 & 4,931/4,931 & 5,653 & 16.0\\
Opt-Verifier & 2,139 & 4,420/4,931 & 81,458 & 120.7\\
\textbf{ORSMT} & \textbf{2,378} & \textbf{4,931/4,931} & 15,249 & 46.1\\
\bottomrule
\end{tabular*}
\end{table}

\begin{table}[htbp]
\caption{Paired gains and losses for protected ORSMT over 3,034 problems and six generators: 18,192 answer-scored and 18,204 model--code-checked outputs.}
\label{tab:nativepairedcounts}
\centering\footnotesize\setlength{\tabcolsep}{3pt}
\begin{tabular*}{\textwidth}{@{\extracolsep{\fill}}llrr@{}}
\toprule
Metric & Comparator & ORSMT only & Comparator only\\
\midrule
Correct answer & Original & 2,378 & 0\\
Correct answer & ReLoop & 1,680 & 131\\
Correct answer & Opt-Verifier & 2,009 & 1,770\\
Model--code match & Original & 8,186 & 0\\
Model--code match & ReLoop & 7,771 & 3\\
Model--code match & Opt-Verifier & 11,442 & 1\\
\bottomrule
\end{tabular*}
\end{table}

\subsection{Transfer to OptiMUS-0.3 and OptiMind}\label{app:optimindaudit}
This experiment keeps each host system's base model and 3,034 assigned problems fixed, and compares
the released pipeline with the same pipeline augmented by ORSMT checking and counterexample feedback.
A correct answer passes the host's released final-answer evaluator. An artifact match is
the separate reference-blind audit introduced in this paper. The OptiMUS-0.3 audit applies the exact model--code
relation to all 3,034 assignments. ORSMT selects 1,095 code-to-mathematics
reconstructions, 533 mathematics-to-code reconstructions and 310 outputs from certified voting,
with shared fallback on 1,096 assignments. Paired accounting adds 825 artifact matches
and 552 correct answers. Outputs satisfying both rise from 630 to 1,199, adding 569 with no losses.

Table~\ref{tab:optimindreplayaudit} reports the full replay over every OptiMind assignment.

\begin{table}[htbp]
\caption{OptiMind integration on 3,034 assigned problems.}
\label{tab:optimindreplayaudit}
\centering\footnotesize\setlength{\tabcolsep}{4pt}
\begin{tabular*}{\textwidth}{@{\extracolsep{\fill}}lrr@{}}
\toprule
Measure & Original & Original + ORSMT \\
\midrule
Assignments & 3,034 & 3,034 \\
Repaired outputs selected & 0 & 2,902 \\
Code-to-mathematics & 0 & 2,837 \\
Certified voting & 0 & 55 \\
Mathematics-to-code-to-mathematics & 0 & 10 \\
Original outputs retained & 3,034 & 132 \\
Artifact matches & 2,254 & 2,937 \\
Artifact match (\%) & 74.3 & 96.8 \\
Correct answers & 1,885 & 1,904 \\
Correct answers (\%) & 62.1 & 62.8 \\
Correct and matched & 1,678 & 1,901 \\
Correct and matched (\%) & 55.3 & 62.7 \\
New artifact matches & 0 & +683 \\
New correct answers & 0 & +19 \\
New correct and matched outputs & 0 & +223 \\
Call overhead (\%) & 0.0 & 29.9 \\
Token overhead (\%) & 0.0 & 32.0 \\
\bottomrule
\end{tabular*}
\end{table}

ORSMT uses code-to-mathematics reconstruction, certified voting and
mathematics-to-code-to-mathematics reconstruction, then retains the original output for the
remaining assignments. Selection uses candidate evidence. For the original outputs, the exact
extractor supports 558 scalar-formulation pairs: 283 are model--code matches and 275 remain
undecided. The remaining model--code pairs are compared independently by deterministic 7B Qwen and
12B Gemma auditors with the same prompt, temperature zero, top-p one and seed zero. The prompt checks
variable domains, bounds, objective, constraint families, index ranges and constants. A deterministic
32B Qwen auditor is called only when the first two disagree. Equivalent or refuted requires agreement
from two of three; otherwise the pair is unresolved. The auditors do not receive benchmark answers.
Each audit response has a 256-token limit and is parsed from the three requested JSON fields.
Selected code-to-mathematics and mathematics-to-code-to-mathematics routes count as artifact matches
only after exact recapture and relation checking; outputs from certified voting retain their ORSMT
evidence.

\begin{evidencecase}{Prompt for the reference-blind OptiMind artifact audit}
\footnotesize\raggedright
\textbf{System:} Audit whether a written optimization model and its final Python solver program are
semantically equivalent. Compare only the supplied written model and final code. Never use or infer
a benchmark answer. Check decision-variable domains and bounds, objective sense and coefficients,
every constraint family including index ranges, strict inequalities and constants. Ignore variable
names and syntactic rewrites when a one-to-one correspondence is clear. Return JSON only with
\texttt{state}, \texttt{confidence} and \texttt{reason}. Use \texttt{EQUIVALENT} only when all model
components match, \texttt{REFUTED} only for a concrete mismatch and \texttt{UNRESOLVED} otherwise.

\smallskip
\textbf{User:} \texttt{WRITTEN\_MODEL:} [assistant modeling text through the final code turn]

\texttt{FINAL\_CODE:} [last Python code block]
\end{evidencecase}

Table~\ref{tab:optimindbysource} partitions correct answers and artifact matches by benchmark and
confirms full population accounting. Paired accounting adds 683 artifact matches, 19 correct
answers and 223 outputs satisfying both with no losses.
\begin{table}[htbp]
\caption{OptiMind host integration counts by benchmark over 3,034 assigned problems.}
\label{tab:optimindbysource}
\centering\footnotesize\setlength{\tabcolsep}{2pt}
\begin{tabular*}{\textwidth}{@{\extracolsep{\fill}}lrrrrr@{}}
\toprule
Benchmark & Problems & \shortstack{Artifact match\\Original} & \shortstack{Artifact match\\Original + ORSMT} &
\shortstack{Correct answer\\Original} & \shortstack{Correct answer\\Original + ORSMT} \\
\midrule
MAMO easy & 652 & 534 & 652 & 556 & 575 \\
Text2Opt template & 630 & 373 & 581 & 300 & 300 \\
OptiBench & 422 & 388 & 421 & 343 & 343 \\
Text2Opt unstructured & 248 & 128 & 231 & 93 & 93 \\
NLP4LP & 242 & 230 & 242 & 175 & 175 \\
NL4Opt & 230 & 222 & 230 & 209 & 209 \\
MAMO complex & 211 & 158 & 208 & 73 & 73 \\
OptMATH & 132 & 56 & 121 & 34 & 34 \\
IndustryOR & 98 & 64 & 93 & 30 & 30 \\
LogiOR & 76 & 39 & 71 & 34 & 34 \\
BWOR & 75 & 47 & 69 & 28 & 28 \\
ComplexOR & 18 & 15 & 18 & 10 & 10 \\
\midrule
Total & 3,034 & 2,254 & 2,937 & 1,885 & 1,904 \\
\bottomrule
\end{tabular*}
\end{table}

\subsection{Source requirements after repair}\label{app:sourcereview}
An artifact match measures agreement between the generated mathematical and executable models;
source requirement checks separately measure whether that shared formulation contains the admitted
requirements from the original problem.
Two separate cohorts evaluate source requirements and revision. The 654-pair audit admits 4,614 typed
obligations and yields 293 complete certificates, 54 explicit-source refutations, 75 symmetry-only
refutations and 232 undecided pairs. Tables~\ref{tab:reviewfunnel}--\ref{tab:rowadoption} follow 1,452 finite repair
candidates and the 231 adopted row revisions, rerunning original and proposed programs.

\subsubsection{Prompts for source review and revision}
Both reviews receive the original problem, the checked model and its computed optimal plan. The plan
review asks for one overall requirement-by-requirement judgment. The row review makes the same task
explicit for every constraint, variable declaration and objective term. A flagged output is revised
from the reported findings and then reviewed again. The computed plan is evidence about the candidate
model; no benchmark answer or reference objective is included.

\begin{evidencecase}{Prompts for plan review, row review and source revision}
\footnotesize\raggedright
\textbf{Common system:} You are an operations-research modeling expert.

\smallskip
\textbf{Inputs to both reviews:} [problem statement] [certified mathematical model] [computed
objective value] [up to 200 nonzero variables in the computed optimal plan]

\smallskip
\textbf{Plan review instruction:} Check the model against the problem statement requirement by requirement:
every data value and coefficient, every bound and integrality condition, every stated constraint,
and the objective and its direction. Use the supplied optimal plan to test the requirements it
touches. Return exactly \texttt{NO\_VIOLATION} when all requirements are represented; otherwise list
each misstated or omitted requirement, quoting its sentence and naming the variables involved. Do
not write code.

\smallskip
\textbf{Row review instruction:} For each constraint, quote the sentence it implements and compare every
coefficient and right-hand side. Check every variable type and bound, list each source requirement
implemented by no row, and compare the objective coefficients and direction. Mark discrepancies as
\texttt{MISMATCH} or \texttt{MISSING}, then end with \texttt{VERDICT: NO\_VIOLATION} or
\texttt{VERDICT: VIOLATION}. Do not write code.

\smallskip
\textbf{Revision inputs:} [problem statement] [previous program] and one of [irreducible infeasible
rows and bounds], [unbounded model], or [certified model, computed plan and review findings].

\textbf{Revision instruction:} The problem statement is the implementation target. Correct the
program so that it implements the problem statement. Return complete corrected code only.
\end{evidencecase}

Both reviews use temperature zero. Plan review allows 1,024 output tokens and row review allows
2,048; revision allows 3,072. The plan field contains at most 200 nonzero variables. Plan review is
flagged unless its reply begins with \texttt{NO\_VIOLATION}. Row review uses the final explicit
\texttt{VERDICT}; if that line is absent, the presence of \texttt{MISMATCH} or \texttt{MISSING}
flags the output. At most 40 flagged lines are passed to revision, with the final 1,500 characters
used as a fallback when no such line is present. These parsing rules are fixed before outcome scoring.

\begin{table}[htbp]
\caption{Source review outcomes on 1,452 finite plans before revision.}
\label{tab:reviewfunnel}
\centering\footnotesize\setlength{\tabcolsep}{3pt}
\begin{tabular*}{\textwidth}{@{\extracolsep{\fill}}lrr@{}}
\toprule
Review outcome & Plan review & Row review \\
\midrule
Initially flagged plans & 20 & 783 \\
Flagged among 576 incorrect answers & 11 & 410 \\
Flagged among 876 correct answers & 9 & 373 \\
Plan revisions adopted after recheck & 8 & 231 \\
\bottomrule
\end{tabular*}
\end{table}

\begin{table}[H]
\caption{Correct answers and execution outcomes on 231 adopted row review revisions.}
\label{tab:rowadoption}
\centering\footnotesize\setlength{\tabcolsep}{3pt}
\begin{tabular*}{\textwidth}{@{\extracolsep{\fill}}lrrr@{}}
\toprule
Outcome & Original & Revised & Net change (outputs) \\
\midrule
Correct answer from captured values & 134 & \textbf{174} & +40 \\
Correct answer without exceptions & 133 & \textbf{171} & +38 \\
Correct answer from explicit output & 88 & \textbf{118} & +30 \\
Program exceptions & 11 & \textbf{4} & $-7$ \\
\bottomrule
\end{tabular*}
\end{table}

\subsection{Residual constraint errors after repair}
\begin{figure}[htbp]
\centering
\includegraphics[width=\textwidth]{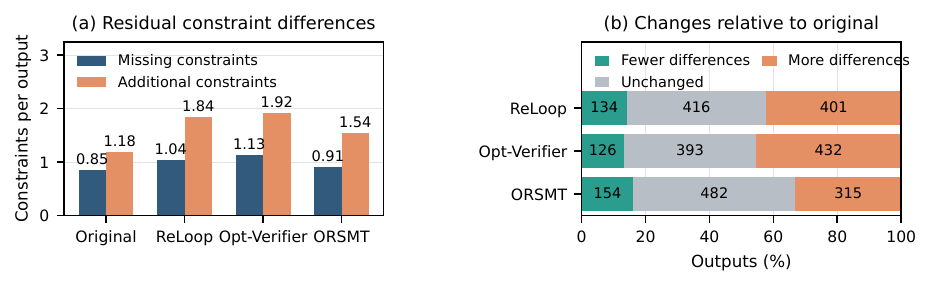}
\caption{Constraint differences on 951 paired outputs. Panel (b) labels give output counts.}
\label{fig:constraintchanges}
\end{figure}

Table~\ref{tab:repairresidual} measures residual missing and additional constraints on the same
fully decided outputs. Figure~\ref{fig:constraintchanges} gives both row totals and the direction
of change per output. All three repair methods introduce constraint differences relative to the original.
ORSMT introduces the fewest and leaves fewer residual differences than either baseline. These are proposals
before protected adoption.
\begin{table}[H]
\caption{Residual constraint differences in revision proposals before protected adoption
on 951 fully decided outputs.}
\label{tab:repairresidual}
\centering\footnotesize\setlength{\tabcolsep}{3pt}
\begin{tabular*}{\textwidth}{@{\extracolsep{\fill}}lrrrrr@{}}
\toprule
Pipeline & \shortstack{Outputs with\\missing\\constraints} & \shortstack{Outputs with\\additional\\constraints} & \shortstack{Outputs without\\constraint\\differences} & \shortstack{Constraint\\precision (\%)} & \shortstack{Constraint\\recall (\%)} \\
\midrule
Original & 415 & 546 & 256 & 72.0 & 78.6 \\
\midrule
ReLoop & 523 & 662 & 176 & 60.6 & 73.6 \\
Opt-Verifier & 544 & 673 & 159 & 58.6 & 71.3 \\
ORSMT & \textbf{462} & \textbf{614} & \textbf{227} & \textbf{66.4} & \textbf{77.0} \\
\bottomrule
\end{tabular*}
\end{table}

The paired row comparison holds the original mathematical artifact fixed. Sequential entailment removes
redundant rows on each side. Missing constraints are mathematical basis rows not entailed by the
program, and additional constraints are program basis rows not entailed by the mathematics.
Constraint recall and precision divide the entailed rows by their respective basis sizes.
The complete comparison requires every relevant row query to decide. Undefined denominators
remain outside conditional precision and recall rates.

The experiment checks 7,620 original and revised pairs with exact rational arithmetic and integer
domains, allowing 200 rows per artifact, five seconds per row and 90 seconds per pair. It uses
12 CPU workers and finishes in 66.0 seconds without language model calls. Eight original records
with unresolved row queries are excluded from complete decisions. The shared intersection contains
951 outputs. Table~\ref{tab:matchingdrift} connects residual differences to apparently correct values.

\begin{table}[htbp]
\caption{Correct answers and changed constraints on the same 951 fully decided outputs.}
\label{tab:matchingdrift}
\centering\footnotesize\setlength{\tabcolsep}{3pt}
\begin{tabular*}{\textwidth}{@{\extracolsep{\fill}}lrrrr@{}}
\toprule
Method & \shortstack{Correct\\answers} & \shortstack{Correct answers with\\constraint differences} & \shortstack{Correct answers with\\more differences} & \shortstack{New differences\\among 256 originals\\without differences} \\
\midrule
ORSMT & 556 & 456 & 215 & 116 \\
ReLoop & 605 & 505 & 260 & 150 \\
Opt-Verifier & 661 & 556 & 303 & 161 \\
\bottomrule
\end{tabular*}
\end{table}

Appendix~\ref{app:workedexamples} expands the formal counterexamples and repair outcomes.
The cases include missed region differences, shared omissions, domain errors, redundant rows
and an incorrect negative solution judgment.

\subsection{Results by benchmark and generator}
Table~\ref{tab:jointdetail1} reports outputs satisfying both conditions for every
benchmark and generator under protected selection. Each cell requires a correct answer and a
model--code match. The problem column gives the denominator
for each generator. These counts resolve the benchmark and generator summaries in the main text
and sum to 7,308 outputs satisfying both without averaging benchmark percentages. Small generators gain strongly on the formulation
benchmarks, while the lower counts on the two Text2Opt families and OptMATH identify benchmarks where fewer outputs satisfy both.
\begin{table}[H]
\caption{Outputs with a correct answer and model--code match by benchmark and generator after protected selection, over 18,192 outputs with numeric references.}
\label{tab:jointdetail1}
\centering\footnotesize\setlength{\tabcolsep}{3pt}
\begin{tabular*}{\textwidth}{@{\extracolsep{\fill}}lrlrrrrrr@{}}
\toprule
Benchmark & Problems & Stage & 1.5B & 3B & 7B & 14B & 32B & 72B \\
\midrule
MAMO easy & 652 & Original & 32 & 81 & 388 & 422 & 580 & 570 \\
 &  & ORSMT & 271 & 396 & 485 & 551 & 596 & 575 \\
Text2Opt template & 630 & Original & 0 & 0 & 6 & 32 & 15 & 79 \\
 &  & ORSMT & 111 & 55 & 45 & 80 & 82 & 110 \\
OptiBench & 422 & Original & 24 & 55 & 142 & 139 & 253 & 277 \\
 &  & ORSMT & 117 & 218 & 236 & 279 & 299 & 290 \\
Text2Opt unstructured & 248 & Original & 0 & 2 & 16 & 61 & 91 & 92 \\
 &  & ORSMT & 15 & 17 & 22 & 64 & 93 & 98 \\
NLP4LP & 242 & Original & 13 & 20 & 64 & 65 & 138 & 143 \\
 &  & ORSMT & 68 & 109 & 125 & 144 & 158 & 156 \\
NL4Opt & 230 & Original & 6 & 20 & 80 & 76 & 150 & 165 \\
 &  & ORSMT & 60 & 141 & 153 & 169 & 181 & 176 \\
MAMO complex & 211 & Original & 0 & 2 & 29 & 29 & 9 & 53 \\
 &  & ORSMT & 16 & 62 & 57 & 44 & 17 & 63 \\
OptMATH & 132 & Original & 1 & 0 & 1 & 2 & 10 & 7 \\
 &  & ORSMT & 8 & 9 & 6 & 8 & 18 & 13 \\
IndustryOR & 98 & Original & 0 & 1 & 10 & 7 & 26 & 20 \\
 &  & ORSMT & 7 & 13 & 18 & 12 & 27 & 22 \\
LogiOR & 76 & Original & 0 & 0 & 1 & 7 & 9 & 11 \\
 &  & ORSMT & 3 & 6 & 4 & 12 & 12 & 14 \\
BWOR & 75 & Original & 1 & 0 & 6 & 5 & 11 & 14 \\
 &  & ORSMT & 4 & 4 & 7 & 9 & 14 & 15 \\
ComplexOR & 18 & Original & 0 & 0 & 1 & 5 & 7 & 6 \\
 &  & ORSMT & 3 & 7 & 5 & 6 & 9 & 9 \\
\bottomrule
\end{tabular*}
\end{table}

\section{Semantic Criterion Evaluation}\label{app:criterionvalidity}
\subsection{Controlled semantic edits and preserving rewrites}
Table~\ref{tab:appcriterion} applies five semantic edit operators and five preserving rewrites.

\begin{table}[htbp]
\caption{Controlled semantic edits and preserving rewrites over 1,360 artifact pairs.}
\label{tab:appcriterion}
\centering\footnotesize\setlength{\tabcolsep}{2pt}
\begin{tabular*}{\textwidth}{@{\extracolsep{\fill}}P{0.27\textwidth}rrrr@{}}
\toprule
Edit operator & Applicable & Proved different & Proved equivalent & Unresolved\\
\midrule
Drop constraint & 144 & 123 & 19 & 2\\
Tighten right hand side & 130 & 121 & 9 & 0\\
Flip relation & 129 & 128 & 0 & 1\\
Relax integrality & 94 & 68 & 20 & 6\\
Perturb coefficient & 143 & 122 & 19 & 2\\
\midrule
All edits & 640 & 562 & 67 & 11\\
\bottomrule
\end{tabular*}
\par\medskip
\begin{tabular*}{\textwidth}{@{\extracolsep{\fill}}lrrr@{}}
\toprule
Preserving rewrite & Pairs & Certified equivalent & Refuted or unresolved\\
\midrule
Scale one constraint & 144 & 144 & 0\\
Duplicate one constraint & 144 & 144 & 0\\
Reorder variables and constraints & 144 & 144 & 0\\
Rename variables & 144 & 144 & 0\\
Positively scale the objective & 144 & 144 & 0\\
\midrule
All preserving rewrites & 720 & 720 & 0\\
\bottomrule
\end{tabular*}
\end{table}

Exact region and preference queries label the resulting model relations. Of 640 applicable edits,
562 are different, 67 are equivalent and 11 remain unresolved. Each of the five preserving operators
contributes 144 equivalent pairs. Outcome value compares the optimized value, mapped optimum compares
optimizing decisions after variable mapping, structural matching compares model graphs, and
falsification samples separating points. Table~\ref{tab:appcriteriondecisions} reports unresolved
outcomes instead of dropping them. These controls evaluate artifact semantics under the formal
reference protocol.

\begin{table}[htbp]
\caption{Criterion decisions over 562 proved different and 787 proved equivalent model pairs.}
\label{tab:appcriteriondecisions}
\centering\footnotesize\setlength{\tabcolsep}{2pt}
\begin{tabular*}{\textwidth}{@{\extracolsep{\fill}}P{0.27\textwidth}rrrrrr@{}}
\toprule
& \multicolumn{3}{c}{Proved different} & \multicolumn{3}{c}{Proved equivalent}\\
\cmidrule(lr){2-4}\cmidrule(lr){5-7}
Criterion & Detected & Missed & Unresolved & Accepted & Falsely flagged & Unresolved\\
\midrule
Outcome value & 313 & 197 & 52 & 613 & 138 & 36\\
Mapped optimum & 358 & 152 & 52 & 569 & 182 & 36\\
Structural graph & 562 & 0 & 0 & 288 & 499 & 0\\
Falsification with one point & 286 & 276 & 0 & 785 & 2 & 0\\
Falsification with four points & 427 & 135 & 0 & 783 & 4 & 0\\
\bottomrule
\end{tabular*}
\end{table}

\begin{table}[htbp]
\caption{Rates and 95\% Wilson intervals on 562 different and 787 equivalent pairs.}
\label{tab:criterionvalidity}
\centering\footnotesize\setlength{\tabcolsep}{3pt}
\begin{tabular*}{\textwidth}{@{\extracolsep{\fill}}lrrr@{}}
\toprule
Criterion & \shortstack{Changed pairs detected (\%)\\with 95\% interval} & \shortstack{Equivalent pairs accepted (\%)\\with 95\% interval} & \shortstack{Pairs\\decided (\%)} \\
\midrule
Outcome value & 55.7 [51.6, 59.7] & 77.9 [74.9, 80.7] & 93.5 \\
Mapped optimum & 63.7 [59.6, 67.6] & 72.3 [69.1, 75.3] & 93.5 \\
Structural graph & 100.0 [99.3, 100.0] & 36.6 [33.3, 40.0] & 100.0 \\
Falsification with one point & 50.9 [46.8, 55.0] & 99.7 [99.1, 99.9] & 100.0 \\
Falsification with four points & 76.0 [72.3, 79.3] & 99.5 [98.7, 99.8] & 100.0 \\
Exact reference relation & 100.0 [99.3, 100.0] & 100.0 [99.5, 100.0] & 100.0 \\
\bottomrule
\end{tabular*}
\end{table}

Four-point falsification gives the best balance between detecting changed pairs and accepting equivalent pairs.
Using the exact semantic relation as the reference partition, we compare weaker evaluation criteria on the proved model pairs.

\subsection{Large edit evaluation}
The larger panel tests which model changes remain invisible to objective values or to a single
region inclusion. Table~\ref{tab:appedit} combines the exact partition and the corresponding miss
rates by operator. Outcome value misses 4,721 proved changes and also falsely separates 13
preserving edits. A test in a single inclusion direction misses 5,738 proved changes. Constraint deletion
and integrality relaxation often preserve the optimum, while tightening requires the opposite
inclusion direction to expose excluded decisions.

\begin{table}[htbp]
\caption{Partitions and miss rates among 12,831 edited pairs. Rates use proved changes in each row.}
\label{tab:appedit}
\centering\footnotesize\setlength{\tabcolsep}{3pt}
\begin{tabular*}{\textwidth}{@{\extracolsep{\fill}}lrrrrrr@{}}
\toprule
Edit operator & Applicable & \shortstack{Proved\\different} & \shortstack{Proved\\equivalent} & Unresolved & \shortstack{Value test\\miss rate (\%)} & \shortstack{One direction\\miss rate (\%)} \\
\midrule
Drop constraint & 2,926 & 2,619 & 298 & 9 & 53.5 & 0.0 \\
Tighten right hand side & 2,489 & 2,330 & 153 & 6 & 36.5 & 100.0 \\
Flip relation & 2,463 & 2,420 & 33 & 10 & 13.2 & 96.6 \\
Relax integrality & 2,028 & 1,240 & 721 & 67 & 63.6 & 0.0 \\
Perturb coefficient & 2,925 & 2,172 & 743 & 10 & 62.7 & 49.3 \\
All edits & 12,831 & 10,781 & 1,948 & 102 & 43.8 & 53.2 \\
\bottomrule
\end{tabular*}
\end{table}

\section{Worked Inconsistency and Repair Cases}\label{app:workedexamples}

The following four recorded benchmark cases connect the original artifacts, verification evidence,
repair and subsequent checks. They connect individual mechanisms to exact variable domains,
objective directions and constraint coefficients.

\subsection{A ReLoop pass with inconsistent feasible regions}
\label{app:industryexample}

\begin{figure}[H]
\begin{evidencecase}{IndustryOR 71: a ReLoop pass with opposite model and code constraints}
\begin{minipage}[t]{0.475\linewidth}
\vspace{0pt}\footnotesize
\textcolor{verdictblue}{\textbf{Mathematical model.}}
The 14 billion parameter generator returns four nonnegative integer quantities
\(\mathrm{boat},\mathrm{plane},\mathrm{train},\mathrm{truck}\).
Both artifacts use the same objective and resource limits.

\textcolor{verdictorange}{\textbf{Written restriction.}}
\[
\mathrm{plane}-\mathrm{boat}\geq0.
\]
\textcolor{verdictorange}{\textbf{Solver code.}}\par
\mbox{\texttt{m.addConstr(boat - plane >= 0)}}

\textcolor{verdictblue}{\textbf{Returned result.}}
The reported objective is 623, equal to the benchmark reference.
ReLoop returns \textsc{Verified} with confidence 1.
\end{minipage}\hfill
\begin{minipage}[t]{0.475\linewidth}
\vspace{0pt}\footnotesize
\textcolor{verdictorange}{\textbf{SMT counterexamples.}}
In the decision order shown at left, \((0,1,0,0)\) satisfies every written row but violates
the executed restriction. The point \((1,0,1,0)\) satisfies every executed row but violates
the written restriction. Both region inclusions are \mbox{\textsc{Refuted}}.

\smallskip
\textcolor{verdictgreen}{\textbf{Mathematics as the target.}}
Regenerate the program with the written restriction.

\smallskip
\textcolor{verdictgreen}{\textbf{Code as the target.}}
Recover the mathematical model with the executed restriction.

\smallskip
\textcolor{verdictgreen}{\textbf{Recheck.}}
Each reconstructed pair has identical objectives and two \mbox{\textsc{Unsat}} region difference queries.
\end{minipage}
\end{evidencecase}
\caption{A ReLoop pass with exact counterexamples and two verified reconstruction directions.}
\label{fig:relooppass}
\end{figure}

IndustryOR problem 71 at 14B returns the reference value 623. The mathematical model and program
each contain four variables and five constraints. ReLoop accepts the program's execution, and the
reported solution passes feasibility, value binding and optimality checks under both models.
The models nevertheless permit different decisions. Figure~\ref{fig:relooppass} shows the original
judgment, the two witnesses and the two repair directions.

\textbf{Complete compared models.} Let \(b,p,r,t\) denote boat, plane, train and truck quantities.
All four are nonnegative integers. Both artifacts maximize
\begin{equation}
5t+10p+8b+7r
\label{eq:industryobjective}
\end{equation}
and contain the four common constraints
\begin{equation}
\begin{aligned}
12t+20p+15b+10r&\leq890,\\
6t+3p+5b+4r&\leq500,\\
t-r&\leq0,\\
b-r&\leq0.
\end{aligned}
\label{eq:industrycommon}
\end{equation}
The remaining mathematical row is \(p-b\geq0\), while the program adds \(b-p\geq0\).
Their objective coefficients, row counts and variable correspondence agree. The different
inequality restricts the feasible region in both directions.

\textbf{Witness checks.} The decision \((b,p,r,t)=(0,1,0,0)\) satisfies every mathematical
constraint. Its wood and steel use are 20 and three, and its executable boat row has value negative
one. Conversely, \((1,0,1,0)\) satisfies the program's complete model, uses 25 units of wood and nine
of steel, and violates the written boat row by one. These decisions are witnesses to the two
set differences, independently of the optimum shared by the models.

\begin{table}[htbp]
\caption{Exact witness evaluations for IndustryOR 71 at 14B.}
\label{tab:industrywitnesses}
\centering\small\setlength{\tabcolsep}{4pt}
\begin{tabular*}{\textwidth}{@{\extracolsep{\fill}}lrr@{}}
\toprule
Quantity or check & Mathematical witness & Executable witness\\
\midrule
Decision \((b,p,r,t)\) & \((0,1,0,0)\) & \((1,0,1,0)\)\\
Wood use, limit 890 & 20 & 25\\
Steel use, limit 500 & 3 & 9\\
Mathematical row \(p-b\), lower bound 0 & 1 & \(-1\)\\
Executable row \(b-p\), lower bound 0 & \(-1\) & 1\\
Objective value & 10 & 15\\
All mathematical constraints satisfied & Yes & No\\
All executable constraints satisfied & No & Yes\\
\bottomrule
\end{tabular*}
\end{table}

\textbf{Repair from the mathematical model.} Keeping the mathematical model as the selected
target replaces the executable inequality with \(p-b\geq0\). The compiler also reproduces the
four common rows, integer declarations and objective. Capturing the new program yields the same
complete model. Both difference queries are unsatisfiable, and the objective expressions are equal.
This direction repairs the inconsistent implementation from the available mathematics.

\textbf{Repair from the program.} Keeping the captured program as the selected target exports
\(b-p\geq0\) into the mathematical artifact and retains every other captured field. The regenerated
program and exported model again have unsatisfiable difference queries and equal objectives.
This direction completes or corrects mathematics from an executable model. The source requirements
determine which selected interpretation should be retained.

\begin{evidencecase}{IndustryOR 71: the changed code and the certified result}
\begin{minipage}[t]{0.475\linewidth}\vspace{0pt}\small
\textcolor{verdictblue}{\textbf{Original program row.}}\par\smallskip
\texttt{m.addConstr(boat - plane >= 0)}\par\medskip
\textcolor{verdictgreen}{\textbf{Mathematical target repair.}}\par\smallskip
\texttt{m.addConstr(plane - boat >= 0)}\par\medskip
The variable declarations, objective and four other rows are retained. The program prints its
computed objective and variable assignments after optimization.
\end{minipage}\hfill
\begin{minipage}[t]{0.475\linewidth}\vspace{0pt}\small
\textcolor{verdictgreen}{\textbf{Certification after construction.}}\par\smallskip
Both region difference queries are unsatisfiable. The objective coefficients and direction match.
A decision with value 623 exists, and the query for a value greater than 623 is unsatisfiable.\par\medskip
\textcolor{verdictblue}{\textbf{Effect of the repair.}}\par\smallskip
The model--code relation changes from refuted to certified even though the answer remains correct.
\end{minipage}
\end{evidencecase}

\textbf{Executed reconstruction.} A fresh replay executes the original program and both
reconstructed Gurobi programs. It reads the instantiated variables, bounds, rows and objective from
each program and compares them with the selected target using exact Z3 queries. Both repaired
programs return \((b,p,r,t)=(0,0,89,0)\) with value 623. The returned point satisfies the complete
model and binds to that value, and no better feasible decision exists. This replay checks the
actual constructed programs in addition to the displayed row calculations.

\subsection{Missing source requirements in two generated artifacts}
\label{app:mamoexample}

MAMO easy problem 20 at 1.5B illustrates two distinct repairs. The first restores agreement between
the generated mathematical model and code. The second restores source requirements that both
artifacts omitted and changes the optimal result.

\textbf{Source and generated models.} The problem chooses integer quantities of four tour
packages, with costs 100, 200, 300 and 400. It gives two capacity limits and two minimum quality
requirements. The generated model uses nonnegative integers \(x_1,x_2,x_3,x_4\), minimizes
\begin{equation}
f(x)=100x_1+200x_2+300x_3+400x_4,
\label{eq:mamocost}
\end{equation}
and contains only \(x_1+x_2\leq50\). The program also includes \(x_3+x_4\leq80\).
Both artifacts omit the quality requirements. Their reported decision is zero in every coordinate,
with objective zero. The recorded report has no explicit solver status.

\begin{evidencecase}{MAMO easy 20: disagreement and missing requirements}
\begin{minipage}[t]{0.475\linewidth}\vspace{0pt}\small
\textcolor{verdictblue}{\textbf{Mathematical artifact.}}\par
\(x_i\in\mathbb Z_{\geq0}\), minimize \(f(x)\),\par
\(x_1+x_2\leq50\).\par\medskip
\textcolor{verdictorange}{\textbf{Additional code row.}}\par
\(x_3+x_4\leq80\).\par\medskip
\textcolor{verdictorange}{\textbf{Exact counterexample.}}\par
\((0,0,81,0)\) satisfies the written model and violates the additional capacity row by one.
\end{minipage}\hfill
\begin{minipage}[t]{0.475\linewidth}\vspace{0pt}\small
\textcolor{verdictblue}{\textbf{Omitted source requirements.}}\par
\(0.6x_1+0.8x_2\geq20\),\par
\(0.7x_3+0.9x_4\geq30\).\par\medskip
\textcolor{verdictorange}{\textbf{Old report.}}\par
\((0,0,0,0)\) has zero cost and is optimal under both incomplete artifacts, while violating both
source quality requirements.
\end{minipage}
\end{evidencecase}

\textbf{Artifact repair.} Adding the second capacity row to the written model makes its region
identical to the program's region. The objective and domains were equal before this change and
remain equal. The zero decision remains optimal. This repair resolves the artifact discrepancy
without changing the interpretation shared by both generated models.

\textbf{Source repair.} The reviewed source formulation contains both capacities and both
quality requirements. Clearing the decimal denominators gives the complete repaired model
\begin{equation}
\begin{aligned}
\min\quad &100x_1+200x_2+300x_3+400x_4\\
\text{subject to}\quad
 &x_1+x_2\leq50, &x_3+x_4&\leq80,\\
 &3x_1+4x_2\geq100, &7x_3+9x_4&\geq300,\\
 &x_1,x_2,x_3,x_4\in\mathbb Z_{\geq0}.&&
\end{aligned}
\label{eq:mamosourcerepair}
\end{equation}
The integer coefficients represent the source decimals exactly. The old zero report fails both
quality checks. The repaired decision \((34,0,43,0)\) has value 16,300 and matches the released
answer. A fresh exact query for a feasible decision with lower cost is unsatisfiable.

\begin{table}[htbp]
\caption{Feasibility and value checks for the repaired MAMO easy 20 decision.}
\label{tab:mamorepairchecks}
\centering\small\setlength{\tabcolsep}{4pt}
\begin{tabular*}{\textwidth}{@{\extracolsep{\fill}}lrr@{}}
\toprule
Constraint or objective & Evaluated value & Required relation\\
\midrule
First capacity, \(x_1+x_2\) & 34 & At most 50\\
Second capacity, \(x_3+x_4\) & 43 & At most 80\\
First quality, \(3x_1+4x_2\) & 102 & At least 100\\
Second quality, \(7x_3+9x_4\) & 301 & At least 300\\
Objective, \(f(x)\) & 16,300 & Equals reported value\\
Feasible decision with smaller objective & None & Query unsatisfiable\\
\bottomrule
\end{tabular*}
\end{table}

\textbf{Result of each intervention.} Model--code verification first restores the missing capacity
row, changing the mathematical/executable row counts from one/two to two/two and the model--code
relation from refuted to certified; the objective remains zero and the answer remains incorrect.
Source review then restores both quality rows in both artifacts, giving four/four rows, objective
16,300 and a correct answer. Code construction propagates the completed mathematical model into the
implementation; the returned decision passes the feasibility and value checks above.

\subsection{A domain mismatch and redundant constraints}
\label{app:complexexample}

ComplexOR problem 3 at 1.5B shows why row agreement alone is insufficient. The written model uses
nonnegative integers \(x,y\), maximizes \(x+y\), and repeats the same row three times:
\begin{equation}
c_1:x+y\leq2,\qquad c_2:x+y\leq2,\qquad c_3:x+y\leq2.
\label{eq:complexrows}
\end{equation}
The code uses the same objective and rows, but declares \(x\) binary. The decision \((2,0)\)
satisfies all rows of both artifacts and violates only the executable domain. The original report
\((0,2)\), with value two, is feasible, correctly valued and optimal under both models.

\begin{evidencecase}{ComplexOR 3: domain repair and a smaller equivalent constraint set}
\begin{minipage}[t]{0.475\linewidth}\vspace{0pt}\small
\textcolor{verdictorange}{\textbf{Domain difference.}}\par
Mathematics: \(x\in\mathbb Z_{\geq0}\).\par
Code: \(x\in\{0,1\}\).\par\medskip
\textcolor{verdictorange}{\textbf{Separating decision.}}\par
At \((x,y)=(2,0)\), each row has left hand side two, while the binary domain rejects \(x=2\).\par\medskip
\textcolor{verdictgreen}{\textbf{Domain repair.}}\par
Declare the code variable as a nonnegative integer and check both inclusions again.
\end{minipage}\hfill
\begin{minipage}[t]{0.475\linewidth}\vspace{0pt}\small
\textcolor{verdictblue}{\textbf{Redundancy evidence.}}\par
For each duplicate, the other two rows entail it. Delete \(c_3\), then test and delete \(c_2\)
against the remaining model.\par\medskip
\textcolor{verdictgreen}{\textbf{Final model.}}\par
One row \(c_1:x+y\leq2\) remains. The integer domains and objective are unchanged. The maximum
value remains two.
\end{minipage}
\end{evidencecase}

\textbf{Checked stages.} The original three-row mathematical model is the target, whereas the
three-row executable model is refuted because it makes \(x\) binary. Restoring \(x\) as a
nonnegative integer certifies the pair. Sequentially removing one and then two duplicate rows
retains certification, leaving one essential row with the same domains and objective.

\begin{figure}[htbp]
\begin{minipage}[t]{0.52\textwidth}\vspace{0pt}
\centering
\includegraphics[width=\linewidth]{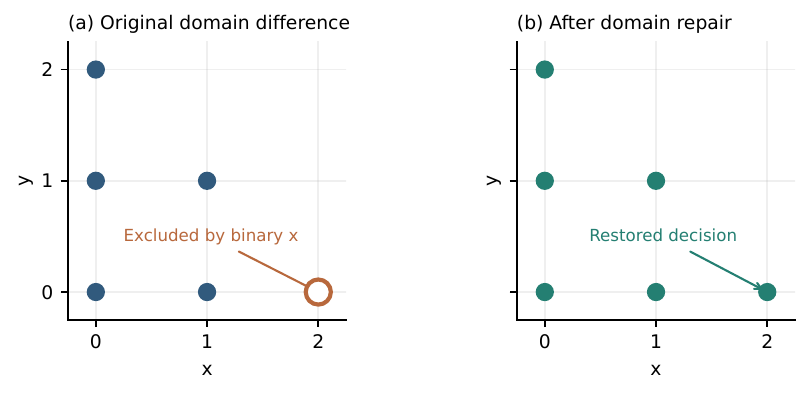}
\caption{ComplexOR 3: the original binary domain excludes the circled integer decision; repair restores it.}
\label{fig:domaincase}
\end{minipage}\hfill
\begin{minipage}[t]{0.44\textwidth}\vspace{0pt}\small
\textbf{Sequential deletion.} To remove \(c_3\), the checker proves that \(c_1,c_2\) and
\(x+y>2\) cannot hold together. After deleting \(c_3\), the same query with \(c_1\) proves
\(c_2\) redundant. The remaining row is essential because \((3,0)\) violates the original bound.
Thus each deletion preserves an equivalent constraint set.

\smallskip
\textbf{Source interpretation.} The source concerns three participants, three cars and an interest
matrix, whereas both artifacts use two variables. Domain and redundancy repair certifies the
generated pair, but restoring the source requires indexed assignment variables and constraints.
The released answer is three; the repaired generated model still has optimum two.
\end{minipage}
\end{figure}

\subsection{An incorrect Opt-Verifier judgment on a certified solution}
\label{app:bworexample}

\begin{evidencecase}{BWOR 18: a correct result receives incorrect revision advice}
\begin{minipage}[t]{0.475\linewidth}\vspace{0pt}\small
\textcolor{verdictgreen}{\textbf{Returned assignment.}}\par
The companies choose stores 3, 2, 1, 4 and 5. The costs are 7, 9, 6, 6 and 6, giving 34 in total.\par\medskip
\textcolor{verdictgreen}{\textbf{Formal checks.}}\par
The model pair is consistent. The decision is feasible, its value is bound correctly, and no
better feasible decision exists in the compared models.
\end{minipage}\hfill
\begin{minipage}[t]{0.475\linewidth}\vspace{0pt}\small
\textcolor{verdictorange}{\textbf{Verifier comment.}}\par
The proposed alternative uses stores 1, 2, 3, 4 and 5. Its listed costs are 4, 9, 12, 6 and 6.
The comment labels this total as 34.\par\medskip
\textcolor{verdictblue}{\textbf{Arithmetic check.}}\par
The listed alternative costs total 37. The alternative supplies no improvement over the certified
value 34.
\end{minipage}
\end{evidencecase}

BWOR problem 18 at 14B assigns five construction companies to five stores. Its recorded model
contains 25 variables and ten constraints. Both region inclusions, objective coefficients and
preferences agree. The reported objective 34 matches the reference value, and the recorded exact
checks certify both region inclusions, equal objective coefficients, equal decision preferences,
feasibility under both models and objective binding. The query for a better feasible decision is
unsatisfiable.

The Opt-Verifier reproduction accepts the structural comparison but returns a negative solution
judgment. Its proposed alternative has store indices \((1,2,3,4,5)\) and costs
\((4,9,12,6,6)\), which sum to 37 rather than the claimed 34. The returned assignment
\((3,2,1,4,5)\) instead costs (7+9+6+6+6=34). The error is in the feedback arithmetic, not
model extraction or variable alignment.

\textbf{Exact decision path and consequence.} The checker proves both region inclusions and equal
objective coefficients, substitutes the returned assignment into every row, binds the reported
value to the evaluated objective and proves that no feasible assignment has value below 34. The
archived comparison judgment uses five language model calls and 44.2 seconds, yet its proposed
counterexample is strictly worse. The exact certificate therefore retains the correct program and
value instead of triggering an unnecessary revision.

\end{document}